\documentstyle[12pt]{amsart}
\newcommand{\OM}{\Omega_f M}
\newcommand{\sm}{\setminus}
\newcommand{\U}{{\Upsilon}}
\newcommand{\E}{{\cal E}}
\newcommand{\J}{{\bf J}}
\newcommand{\Z}{{\Bbb Z}}
\newcommand{\R}{{\Bbb R}}
\newcommand{\li}{\langle i \rangle}
\newcommand{\ti}{\langle \tilde i \rangle}
\author{V.~A.~Vassiliev}

\address{Steklov Math. Inst., Gubkina st.~8, 117966 Moscow, RUSSIA}

\curraddr{MSRI, 1000 Centennial Drive, Berkeley, CA 94720, USA}

\email{vassil@@vassil.mccme.rssi.ru}

\title{On invariants and homology of spaces 
of knots in arbitrary manifolds}
\thanks{Supported in part by  
NSF grant DMS-9022140, 
RFBR (project \# 95-01-00846a) and 
INTAS (project \# 4373). 
}

\begin{document}

\begin{abstract}
The construction of finite-order knot invariants 
in ${\Bbb R}^3,$ based on resolutions of discriminant sets,
can be carried out immediately
to the case of knots in an arbitrary three-dimensional manifold $M$ (may be 
non-orientable)
and, moreover, to the calculation of cohomology groups of spaces
of knots in arbitrary manifolds of dimension $\ge 3$.

Obstructions to the integrability of admissible weight systems
to well-defined knot invariants in $M$ are identified as 
1-dimensional cohomology classes of 
certain generalized loop spaces of $M$. Unlike the case 
$M = \R^3$, these obstructions
can be non-trivial and provide invariants of the
manifold $M$ itself.

The corresponding algebraic machinery allows us to obtain on the level of
the ``abstract nonsense''
some of results and problems of the theory, and 
to extract from others the essential topological 
(in particular, low-dimensional) part.
\end{abstract}

\maketitle

\section*{Introduction}

Finite-order invariants of knots in $\R^3$ appeared in 
\cite{v90} from a topological study of the discriminant subset
of the space of curves in $\R^3$ (i.e., the set of all singular curves).

In \cite{lin}, \cite{kalf}, finite-order invariants of
knots in 3-manifolds satisfying some conditions were considered:
in \cite{lin} it was done for manifolds with $\pi_1=\pi_2=0$,
and in \cite{kalf} for closed irreducible orientable manifolds. In particular,
in \cite{lin} it is shown that the theory of finite-order invariants of knots
in two-connected manifolds is isomorphic to that for $\R^3$;
the main theorem of \cite{kalf} asserts that for every 
closed oriented irreducible 3-manifold $M$
and any connected component of the space $C^\infty(S^1,M)$ the
group of finite-order invariants of knots from this component
contains a subgroup isomorphic to the group of 
finite-order invariants in $\R^3$.
The starting point of these generalizations was the characteristic
property of finite-order invariants in $\R^3$, considered in 
\cite{v90}, \cite {BL}, \cite{BN}: 
the triviality of  certain ``finite differences'' 
of values of these invariants defined with the  
help of the orientation of the ambient manifold, see \S \ 1.2 below.

We show that all the theory of  
finite-order invariants and of the cohomology
of spaces of knots in $\R^3$, based on the study of 
discriminants, can be extended
almost immediately to the case of arbitrary manifolds of dimension
$n \ge 3$ (including non-orientable ones),
although the answers generally are not so easy.

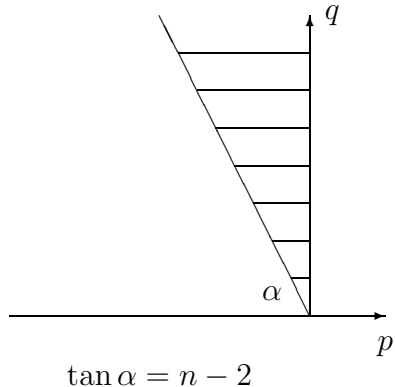
\begin{figure}
\unitlength=1.00mm
\special{em:linewidth 0.4pt}
\linethickness{0.4pt}
\begin{picture}(101.00,50.00)
\put(10.00,10.00){\vector(1,0){50.00}}
\put(50.00,10.00){\vector(0,1){40.00}}
\put(50.00,10.00){\line(-1,2){20.00}}
\put(50.00,15.00){\line(-1,0){2.50}} 
\put(50.00,20.00){\line(-1,0){5.00}}
\put(50.00,25.00){\line(-1,0){7.50}}
\put(50.00,30.00){\line(-1,0){10.00}}
\put(50.00,35.00){\line(-1,0){12.50}}
\put(50.00,40.00){\line(-1,0){15.00}}
\put(50.00,45.00){\line(-1,0){17.50}}
\put(60.00,6.00){\makebox(0,0)[cc]{$p$}}
\put(53.00,50.00){\makebox(0,0)[cc]{$q$}}
\put(45.00,13.00){\makebox(0,0)[cc]{$\alpha$}}
\put(30,2){\makebox(0,0)[cc]{$\tan \alpha = n-2$}}
\end{picture}
\caption{Support of the spectral sequence}
\label{support}
\end{figure}

The cohomology classes of spaces of knots in $M^n$ 
come from a spectral sequence
with support in the wedge
$\{(p,q)| p<0, p(n-2)+q \ge 0\}$, see fig.~\ref{support}. 

In particular, for $n=3$ the 
knot invariants (= 0-dimensional cohomology classes
of such a space) are counted by terms $E_\infty^{-i,i}$ of this sequence.
In this case all elements of terms $E_0^{-i,i}$ 
of our spectral sequence, and some elements of $E_0^{-i,i+1}$ appeared
in \cite{lin} and \cite{kalf} under the names ``singular
knot invariants'' and ``local integrability conditions'';
as in \cite{v90}, the calculation of further terms of the 
spectral sequence is nothing
but the check whether these singular invariants satisfying 
these local conditions can be extended  
to less complicated singular knots or not.

In the case of an arbitrary manifold $M^3$, groups
$E_r^{-j,j+1}$ contain some additional elements, the
``global integrability conditions'', which can provide non-trivial homological
obstructions to this extension, see \S \ 1.5. Our spectral
sequence allows us to write these obstructions explicitly and,
moreover, to be sure that if for some initial data they vanish,
then these data can be extended to a well-defined knot invariant.
Also the use of spectral sequences allows us to prove the main theorem
of \cite{lin} and similar more general comparison theorems
just by the methods of ``abstract nonsense'', and to avoid the
technical difficulties overcome in \cite{lin}, \cite{kalf} 
by the methods from the Yablokova work.
\medskip

The principles of the paper can be formulated in the following 
five statements.
\medskip

{\bf 1. The ``manifold'' part of invariants also is interesting.}

For $n=3$ the cohomology classes coming from the ``principal diagonal''
$\{E_\infty^{-i,i}\}$ of
the spectral 
sequence are, generally speaking, the invariants not of 
knots in $M^3$ but of both the knots and the manifold
$M^3$; they are exactly the ``singular knot
invariants'' in terminology of \cite{lin}, \cite{kalf}. 
As was pointed out in \cite{lin}, they
take values on pairs of the 
form \{a knot in $M$; a path in the 
space of continuous maps $S^1 \to M^3$ connecting this knot to 
a distinguished knot in its homotopy class (and considered up to homotopy)\}.

More generally, for $M$ of arbitrary dimension $n$ let $\OM$ be the space
of smooth maps $S^1 \to M$ and $\Sigma \subset \OM$ the set of all
maps having selfintersections or singular points, so that knot invariants
are elements of the group $H^0(\OM \setminus \Sigma)$.
Then our spectral sequence converges to a subgroup of the relative
cohomology group $H^*(\OM, \OM \setminus \Sigma)$
(if $n>3$ then to all this group:
\begin{equation}
\label{ss}
E_r^{p,q} \to H^{p+q+1}(\OM, \OM \setminus \Sigma)).
\end{equation}

In particular, for $n=3$ elements of the limit group $E_\infty^{-i,i}$ 
of this sequence define 1-dimensional cohomology classes of
$\OM$ (and these classes can be nontrivial);
and elements, defining zero cohomology classes, can be lifted to well-defined
 knot invariants.

Restrictions on manifolds, required in \cite{lin} and \cite{kalf}, 
essentially describe some situations, when all or some of these 
1-cohomology classes are trivial. In these cases the
invariants, provided by the above 
spectral sequence, coincide with these from 
\cite{lin}, \cite{kalf}.

However these cohomology classes are an interesting characteristic
of the manifold $M$ and probably should not be considered separately from
knot invariants. Moreover, the spectral sequence itself (especially
its higher differentials) is a strong invariant of $M$.
By analogy with the main result of \cite{BL}, I wonder whether the 
Jones--Witten--Reshetikhin--Turaev invariants of $M$ can 
be derived from these ones.
\medskip

{\bf 2. One can calculate also higher-dimensional cohomology of spaces
of knots, in particular in manifolds of higher dimensions.}

The simplest such class is the 1-dimensional cohomology class of 
order 1 (i.e. coming from the cell $E^{-1,2}$) of the space
of knots in $\R^3$. It takes nontrivial value on a 1-cycle in the
space of unknots in $\R^3$, in particular proves that this space
is not simple-connected, see \S \ 1.8. 

On the other hand, for knots in $M^n$, $n>3,$ we have 
no problems with the convergence of 
the spectral sequence, see (\ref{ss}).
An attractive (for me) problem is to study explicitly these 
spectral sequences for 
simple-connected 4-dimensional manifolds.
\medskip

{\bf 3. There is an essential theory of finite-order invariants
in non-ori\-en\-tab\-le manifolds.}

The orientability of $M^3$, used in \cite{kalf} and
\cite{lin} for the transversal orientation of the discriminant set,
is unneseccary. First of all, the entire theory can be carried out
without changes to the case of cohomologies and spectral sequences
with coefficients in $\Z_2$, when such a coorientation is useless.
Moreover, as we shall see in \S \ 1.6, even for non-orientable manifolds
our construction can give nontrivial invariants with integer
(or real) values. In this case the basic ``axiomatic'' definition of
order $i$ invariants should be replaced by a more general one, as well as 
the notion of the index of order $i$, which any invariant defines at
a singular knot with $i$ self-intersections.
\medskip

{\bf 4. Resolutions of singularities of the discriminant are useful.}

To show that the intersection (or linking) number with a subvariety
is well defined, one usually tries to prove that this variety 
is regular (up to a set of codimension 2) and
(co)orientable, or (if this is wrong) that it is possible to 
orient all its smooth pieces in such a way that these 
orientations will be compatible
close to singular points, cf. \cite{arnold1}, 
\cite{vll}, \cite{arnold}, \cite{lin}.

Often such considerations can be simplified very much.
Indeed, it is sufficient to show that our variety is the image of a regular
orientable manifold under a proper map, and to define our indices
in terms of the direct image of its fundamental cycle. 
The representation of varieties as such images is provided by
the techniques of singularity resolutions.  If our variety
is a stratum of the discriminant set 
(say, the closure a class of singular knots), then 
there is a ``tautological'' construction of such resolutions.
\medskip

{\bf 5. Spectral sequences also are useful.}

Many comparison theorems of topology and algebraic geometry can be proved in
the following standard way. Instead of comparing homology
groups of two objects, we compare spectral sequences converging to
these groups. The convergence process (and the resulting groups)
can be very complicated, but we do not need to consider it:
it is easier to prove that 
the initial terms of these spectral sequences are naturally isomorphic,
and isomorphisms of
their final terms and these groups will follow automatically. 
Similar considerations
often prove that one of these groups is ``greater'' than the other.

For instance, our spectral sequences are functorial with respect
to the inclusion of manifolds, and hence the coincidence theorem from
\cite{lin} is an immediate corollary of the 
(very easy) comparison of initial terms of corresponding 
spectral sequences for $M^3$ and $\R^3$, induced by any embedding
$\R^3 \to M^3$, see \cite{congr} and \S \ 1.7 below. These considerations 
allow to prove immediately some existence theorems for invariants,
or at least to reduce problems of this kind 
to essential low-dimensional
problems concerning initial terms $E_1$ of these sequences 
(which are more standard than the study of the ``integration process'').

Moreover, our sequences are evidently functorial with respect
to any submersions (in particular coverings) of manifolds of the
same dimension. Another attractive problem: which part of this
functoriality survives for arbitrary smooth maps.
\bigskip

In  section 1 we outline main
features of the spectral sequence and describe in
elementary terms its applications to 
knot invariants; in section 2 we give 
the exact constructions and technical details.
The ``invariant--theoretical'' part of section 1 can be
reformulated and proved in terms not referring to the theory of
spectral sequences,
however we preserve there the standard notation of this theory.  

Also in \S \ 1 I represent four persons: a non-integrable singular knot
invariant (see \S~1.4), an invariant of singular knots with $\ge 2$
crossings, which cannot be extended to invariants of knots with one
crossing (\S \ 1.5), a non-trivial integer first order invariant in
a non-orientable manifold (\S \ 1.6.1), and an invariant of order one,
proving the nontriviality
of the Whitehead link (\S~1.3).
\medskip

I thank very much Efstratia Kalfagianni and Xiao-Song Lin for 
sen\-ding me pre\-prints of their works, and Yasha Eliashberg,
Sergei Matveev and Kolya Mishachev for
related conversations.

Especially I thank MSRI, where the main part of the work was
done, for hospitality and excellent conditions.

\section{Elementary theory and main results}

\subsection{How to overcome the infinitedimensionality.}

Let $M$ be a smooth $n$-di\-men\-si\-onal manifold, $n \ge 3$, $\OM$
the space of all smooth maps $S^1 \to M$, and $\Sigma \subset \OM$
the set of maps having self-intersections or singular points,
so that knots are the elements of $\OM \setminus \Sigma$.

As in \cite{v90}, \cite{book}, we use a sort of the Alexander duality in the 
space $\OM$.
To justify this duality
in the infinitedimensional space, we need to consider a family of 
finitedimensional approximations to this space.
For this, let us embed $M$ into some space $\R^N$ as a regular 
(may be not closed) submanifold; let
$U$ be some open tubular neighborhood of $M$ in $\R^N$, and $\tau: U \to M$  
the corresponding $C^\infty$-smooth projection with open discs
for the fibers. 

For any finitedimensional affine subspace $\Gamma$ of the space of smooth
maps $S^1 \to \R^N$ denote by $\Gamma_U$ its subset consisting of 
maps, whose images belong to $U$. 
For the approximating subsets of the loop
space $\OM$ we will use the sets of maps of the form $\tau \circ f$, 
$f \in \Gamma_U$. 
\medskip

For any such space $\Gamma$ we construct in \S \ 2 
a homological spectral sequence
$$E^r_{p,q}(\Gamma) \to \bar H_{p+q}(\Sigma \cap \Gamma_U)$$
(where
$\bar H_*$ denotes the Borel--Moore homology, i.e. the homology of the
one-point compactification reduced modulo the added point)
in exactly the same way as it was done in \cite{v90}
in the special case $M = \R^3 = \R^N = U$.
Using the formal change of indices
\begin{equation}
\label{inver}
E_r^{p,q} \equiv E^r_{-p,\dim \Gamma -q-1}
\end{equation}
and the Poincar\'e--Lefschetz duality 
$$\bar H_j(\Sigma \cap \Gamma_U) \simeq H^{\dim \Gamma - j}(\Gamma_U, 
\Gamma_U \setminus \Sigma),$$ we convert this sequence to a
cohomological spectral sequence
$$E_r^{p,q}(\Gamma) \to H^{p+q+1}(\Gamma_U, \Gamma_U \setminus \Sigma).$$

If $\Gamma$ satisfies some genericity conditions (see \S \ 2.2), 
then the support
of its term $E_1$ (and hence of all subsequent terms) belongs to the wedge
shown in fig. 1.

\subsection{Stabilization of spectral sequences.}

For any natural $m$, if $\Gamma \subset \Gamma'$ are two subspaces
satisfying these genericity conditions, and the dimension 
of $\Gamma$ is sufficiently large
with respect to $m$, then for any $p \in [-m,0]$ and any $q$
the natural homomorphism $E^{p,q}_r(\Gamma') \to E^{p,q}_r(\Gamma),$
$r \ge 1$,
is well defined, see \S \ 2.3 below. This homomorphism is compatible with all
subsequent differentials $d_r: E_r^{p,q} \to E_r^{p+r, q-r+1}$ 
of the spectral sequence, and its limit version (for $r=\infty$)
is compatible with the map of cohomology groups 
$H^*(\Gamma'_U, \Gamma'_U \setminus \Sigma) \to 
H^*(\Gamma_U,\Gamma_U \setminus \Sigma)$
induced by the identical embedding.
Thus the limit spectral sequence $$E_r^{p,q} \equiv 
\mbox{lim  ind }  E_r^{p,q}(\Gamma)$$
is well defined. 
For $n\ge 4$ this spectral sequence converges exactly to the group 
$H^*(\OM, \OM \setminus \Sigma)$. For $n=3$ the similar statement is
not proved (and probably is wrong, at least for sufficiently complicated
$M$) because the sequence has infinitely many nontrivial terms
$E_1^{p,q}$ on any line $\{p+q=const \ge 0\}$.

For the most popular case, when $n=3$, $p+q=0$ and $M$ is orientable, 
the group $E_1^{p,q}$ is described in subsection 1.3. The group of finite-order
knot invariants, to which these groups $E_r^{-i,i}$
stabilize, can be characterized in the following standard way.
\medskip

\begin{figure}
\unitlength=1.00mm
\special{em:linewidth 0.4pt}
\linethickness{0.4pt}
\begin{picture}(80,20)
\put(0,6){\vector(1,1){12}}
\put(12,6){\line(-1,1){5.5}}
\put(5.5,12.5){\vector(-1,1){5.5}}
\put(25,6){\vector(1,1){12}}
\put(31,12){\circle*{1.33}}
\put(37,6){\vector(-1,1){12}}
\put(50,6){\line(1,1){5.5}}
\put(56.5,12.5){\vector(1,1){5.5}}
\put(62,6){\vector(-1,1){12}}
\put(6,1){\makebox(0,0)[cc]{$a$}}
\put(6,6){\makebox(0,0)[cc]{$+$}}
\put(31,1){\makebox(0,0)[cc]{$b$}}
\put(56,1){\makebox(0,0)[cc]{$c$}}
\put(56,6){\makebox(0,0)[cc]{$-$}}
\put(23,12){\vector(-1,0){9}}
\put(39,12){\vector(1,0){9}}
\end{picture}
\caption{Local resolutions of a selfintersection point}
\label{cp}
\end{figure}
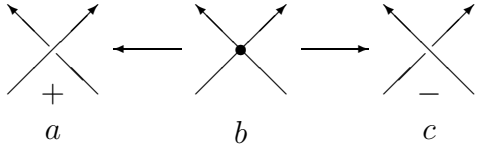

Consider any immersion $\phi: S^1 \to M$,
having a transverse selfintersection point $\phi(x)=\phi(y)$,
$x \ne y \in S^1$. This selfintersection
can be removed in two locally different ways, see fig. \ref{cp}.
Using the orientation of $M$ we can call one of these local perturbations 
{\it positive} and the other {\it negative}, see \S~4.3.1 
in \cite{v90}; opposite
orientations of $M$ define different signs of local perturbations. 
The index of a knot invariant
at the singular knot $\phi$ is defined as its value at the positive 
perturbation minus that at the negative one.

Similarly, if $\phi$ has exactly $j$ different transverse selfintersection
points, then there are $2^j$ different simultaneous perturbations 
of all these points, moving our immersion $\phi$ to non-singular
knots. Any knot invariant associates the {\it index of order 
$j$} to such a singular knot $\phi$: it is equal to the sum of values of the 
invariant over all perturbations, the number of negative local moves
in which is even, minus the sum of similar perturbations with an odd number of
negative local moves.
\medskip

{\sc Definition-Proposition.}
An invariant of knots in an orientable $M^3$ is {\it of order $i$} if 
any of following equivalent conditions is satisfied:

1) it has filtration $\le i$ in our spectral sequence
(i.e., comes from some of its cells $E_\infty^{-l,l}$ with $l \le i$);

2) (see \cite{lin}, \cite{kalf}) all its 
indices of orders $j \ge i$ at all immersions  with $j$ 
selfintersections are equal to $0$.
\medskip

The fact, that these two conditions are equivalent, follows 
immediately from the construction of the spectral sequence.
In \cite{v90}, \cite{book} the first of similar two conditions was
used to define finite-order invariants in $\R^3$, and the second
was mentioned as its ``geometrical interpretation''.
In all the subsequent publications this second condition appears
as the main definition.
\medskip

{\sc Definition.} The local surgery of an immersion $S^1 \to M^3,$ 
connecting two pictures fig. \ref{cp}a and \ref{cp}c,  
is called {\it positive} (respectively, {\it negative}), if it 
replaces the negative resolution of the corresponding 
singular immersion by the positive one (respectively,
positive by negative).

\subsection{Description of the term $E_1^{-i,i}$ of the 
stable spectral sequence for $n=3$.}

{\sc Definition} (see \cite{v90}).
The $[i]$-{\it configuration} (or chord diagram, see \cite{BN}) is any
collection of $2i$ distinct points in $S^1$ partitioned into $i$ pairs.
An $\li$-{\it configuration} is 
any collection of $2i-1$ distinct points in $S^1$ partitioned into $i-2$ pairs
and one triple.
An $i^*$-configuration is any collection of $2i-1$ points in $S^1$
partitioned into $i-1$ pair and one distinguished point $*$.
An $\ti$-configuration is an $\li$-configuration,
in which one point of the triple
is distinguished.

For any symbol $\Upsilon = [i], \li, \ti$ or $i^*$, two
$\Upsilon$-configurations are {\it equivalent} if they can be transformed
into one another by an orientation-preserving homeomorphism of $S^1$.
A map $\phi: S^1 \to M$ {\it respects} the $\Upsilon$-configuration
if it maps all the points of any pair or triple to one point
in $M$ and (in the case of $i^*$-configurations) has zero derivative
at the point $*$.

An $\Upsilon_M$-{\it route} is any pair of the form \{an 
equivalence class of $\Upsilon$-configurations in $S^1$; 
a homotopy class of maps $S^1 \to M$ respecting  
configurations of this class\}.
\medskip

We describe the group $E_1^{-i,i}$ in  two equivalent ways. The first
of them (following \cite{v90}) reflects better the structure of the 
resolution space and can be generalized to the calculation of 
higher-dimensional homology classes, see \S \ 2; the second is more
standard (and is formulated in terms of 4-term relations etc.) 

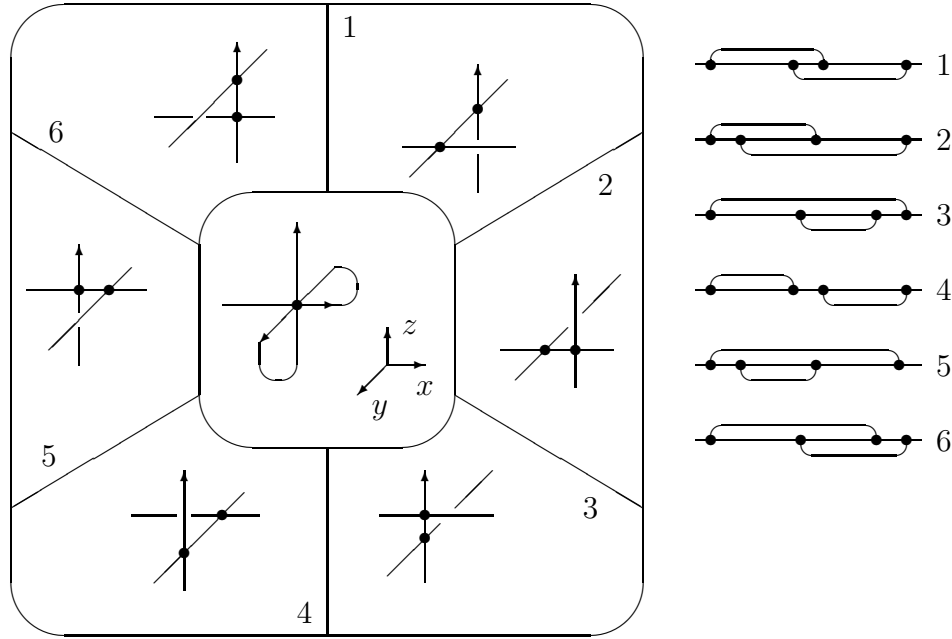
\begin{figure}
\unitlength=1.00mm
\special{em:linewidth 0.4pt}
\linethickness{0.4pt}
\begin{picture}(124.00,85.00)
\put(50.00,37.00){\vector(1,0){5.00}}
\put(50.00,37.00){\vector(0,1){5.00}}
\put(50.00,37.00){\vector(-1,-1){4.00}}
\put(55.00,34.00){\makebox(0,0)[cc]{$x$}}
\put(53.00,42.00){\makebox(0,0)[cc]{$z$}}
\put(49.00,31.00){\makebox(0,0)[cc]{$y$}}
\put(28.00,45.00){\vector(1,0){15.00}}
\put(38.00,40.00){\vector(0,1){16.00}}
\put(43.00,50.00){\vector(-1,-1){10.00}}
\put(38.00,45.00){\circle*{1.33}}
\put(23.00,12.00){\circle*{1.33}}
\put(28.00,17.00){\circle*{1.33}}
\put(55.00,17.00){\circle*{1.33}}
\put(55.00,14.00){\circle*{1.33}}
\put(75.00,39.00){\circle*{1.33}}
\put(71.00,39.00){\circle*{1.33}}
\put(62.00,71.00){\circle*{1.33}}
\put(57.00,66.00){\circle*{1.33}}
\put(30.00,75.00){\circle*{1.33}}
\put(30.00,70.00){\circle*{1.33}}
\put(9.00,47.00){\circle*{1.33}}
\put(13.00,47.00){\circle*{1.33}}
\put(77.00,18.00){\makebox(0,0)[cc]{$3$}}
\put(39.00,4.00){\makebox(0,0)[cc]{$4$}}
\put(5.00,25.00){\makebox(0,0)[cc]{$5$}}
\put(6.00,68.00){\makebox(0,0)[cc]{$6$}}
\put(45.00,82.00){\makebox(0,0)[cc]{$1$}}
\put(79.00,61.00){\makebox(0,0)[cc]{$2$}}
\put(93.00,27.00){\circle*{1.33}}
\put(105.00,27.00){\circle*{1.33}}
\put(119.00,27.00){\circle*{1.33}}
\put(115.00,27.00){\circle*{1.33}}
\put(104.00,27.00){\oval(22.00,4.00)[t]}
\put(112.00,27.00){\oval(14.00,4.00)[b]}
\put(118.00,37.00){\circle*{1.33}}
\put(107.00,37.00){\circle*{1.33}}
\put(93.00,37.00){\circle*{1.33}}
\put(97.00,37.00){\circle*{1.33}}
\put(105.50,37.00){\oval(25.00,4.00)[t]}
\put(102.00,37.00){\oval(10.00,4.00)[b]}
\put(93.00,47.00){\circle*{1.33}}
\put(104.00,47.00){\circle*{1.33}}
\put(108.00,47.00){\circle*{1.33}}
\put(119.00,47.00){\circle*{1.33}}
\put(98.50,47.00){\oval(11.00,4.00)[t]}
\put(113.50,47.00){\oval(11.00,4.00)[b]}
\put(119.00,57.00){\circle*{1.33}}
\put(115.00,57.00){\circle*{1.33}}
\put(105.00,57.00){\circle*{1.33}}
\put(93.00,57.00){\circle*{1.33}}
\put(106.00,57.00){\oval(26.00,4.00)[t]}
\put(110.00,57.00){\oval(10.00,4.00)[b]}
\put(119.00,67.00){\circle*{1.33}}
\put(107.00,67.00){\circle*{1.33}}
\put(93.00,67.00){\circle*{1.33}}
\put(97.00,67.00){\circle*{1.33}}
\put(100.00,67.00){\oval(14.00,4.00)[t]}
\put(108.00,67.00){\oval(22.00,4.00)[b]}
\put(119.00,77.00){\circle*{1.33}}
\put(108.00,77.00){\circle*{1.33}}
\put(104.00,77.00){\circle*{1.33}}
\put(93.00,77.00){\circle*{1.33}}
\put(100.50,77.00){\oval(15.00,4.00)[t]}
\put(111.50,77.00){\oval(15.00,4.00)[b]}
\put(124.00,77.00){\makebox(0,0)[cc]{$1$}}
\put(124.00,67.00){\makebox(0,0)[cc]{$2$}}
\put(124.00,57.00){\makebox(0,0)[cc]{$3$}}
\put(124.00,47.00){\makebox(0,0)[cc]{$4$}}
\put(124.00,37.00){\makebox(0,0)[cc]{$5$}}
\put(124.00,27.00){\makebox(0,0)[cc]{$6$}}
\put(23.00,7.00){\vector(0,1){16.00}}
\put(55.00,8.00){\vector(0,1){15.00}}
\put(75.00,34.00){\vector(0,1){15.00}}
\put(62.00,67.00){\vector(0,1){10.00}}
\put(30.00,64.00){\vector(0,1){16.00}}
\put(9.00,44.00){\vector(0,1){9.00}}
\put(42.00,43.00){\oval(84.00,84.00)[]}
\put(42.00,43.00){\oval(34.00,34.00)[]}
\put(59.00,53.00){\line(5,3){25.00}}
\put(59.00,33.00){\line(5,-3){25.00}}
\put(25.00,33.00){\line(-5,-3){25.00}}
\put(25.00,53.00){\line(-5,3){25.00}}
\put(2.00,47.00){\line(1,0){16.00}}
\put(16.00,50.00){\line(-1,-1){11.00}}
\put(9.00,37.00){\line(0,1){5.00}}
\put(16.00,17.00){\line(1,0){6.00}}
\put(24.00,17.00){\line(1,0){9.00}}
\put(31.00,20.00){\line(-1,-1){12.00}}
\put(50.00,9.00){\line(1,1){7.00}}
\put(59.00,18.00){\line(1,1){5.00}}
\put(64.00,17.00){\line(-1,0){15.00}}
\put(65.00,39.00){\line(1,0){15.00}}
\put(67.00,35.00){\line(1,1){7.00}}
\put(76.00,44.00){\line(1,1){5.00}}
\put(66.00,75.00){\line(-1,-1){13.00}}
\put(52.00,66.00){\line(1,0){15.00}}
\put(62.00,77.00){\line(0,-1){10.00}}
\put(62.00,65.00){\line(0,-1){5.00}}
\put(34.00,79.00){\line(-1,-1){13.00}}
\put(24.00,70.00){\line(-1,0){5.00}}
\put(26.00,70.00){\line(1,0){9.00}}
\put(121.00,27.00){\line(-1,0){30.00}}
\put(91.00,37.00){\line(1,0){30.00}}
\put(121.00,47.00){\line(-1,0){30.00}}
\put(91.00,57.00){\line(1,0){30.00}}
\put(121.00,77.00){\line(-1,0){30.00}}
\put(91.00,67.00){\line(1,0){30.00}}
\put(42.00,1.00){\line(0,1){25.00}}
\put(42.00,60.00){\line(0,1){25.00}}
\put(35.50,40.00){\oval(5.00,10.00)[b]}
\put(43.00,47.50){\oval(6.00,5.00)[r]}
\end{picture}
\caption{Splittings of a triple point}
\label{triple}
\end{figure}

\medskip
{\sc First description.}
 For any natural $i$ the term $E_1^{-i,i}$ is the kernel of a certain operator
$d_0: E_0^{-i,i} \to \tilde E_0^{-i,i+1}$; let us define the elements 
of this operator. Suppose first that the manifold $M$
is orientable or that the coefficient group $G$ is the field $\Z_2$.
Then the group $E^{-i,i}_0$ (respectively, $\tilde E^{-i,i+1}_0$)
is the space of $G$-valued functions on the set of 
$[i]_M$- and $\li_M$-routes (respectively, on the set
of $\ti_M$- and $i^*_M$-routes).

The boundary $d_0(\alpha)$ of the 
generator $\alpha$, corresponding to an $[i]_M$-route, is the sum of 
at most $2i$ generators of the group $\tilde E_0^{-i,i+1}$, corresponding 
to some of segments, into which this configuration divides the circle.
Namely, among these segments there can be ``suspicious'' ones, which
are bounded by the points of one pair of our chord diagram.
To such a segment there corresponds a generator in $d_0(\alpha)$
if and only if the loop in $M$, formed by the image of this segment,
is homotopically trivial; this generator is spanned by the
$i^*_M$-route, represented by the map $S^1 \to M,$ coinciding
with one representing $\alpha$ outside a small neighborhood of this
segment and
replacing this loop by a cusp at its corner point (which will be
the $*$-point of the $i^*_M$-route):
\begin{picture}(46,15)
\put(8,0){\oval(16,16)[lt]}
\put(8,12){\oval(16,16)[lb]}
\put(8,6){\oval(6,4)[r]}
\put(25,6){\makebox(0,0)[cc]{$\to$}}
\put(42,0){\oval(12,12)[lt]}
\put(42,12){\oval(12,12)[lb]}
\end{picture}
. 
To any non-suspicious segment of our
chord diagram, there corresponds the generator in $\tilde E_0^{-i,i+1}$,
equal to the $\ti_M$-route, obtained from $\alpha$
by a degeneration, contracting this segment to a point (which will be
the distinguished point of the triple of the configuration, cf. \cite{v90}).
For instance, the curve with a triple point, shown in the center 
of fig. \ref{triple}, is obtained by such a degeneration from any of 6
singular knots around it.
(However, there is a subtlety here.
Suppose that our
$[i]_M$-route has self-equivalences, i.e. any representing it 
$[i]$-configuration can be transposed into itself by a homeomorphism
of $S^1$ inducing a non-trivial cyclic permutation of its $2i$
vertices in such a way that the composition with this homeomorphism
preserves the homotopy class of maps from our $[i]_M$-route.
Of course, if we contract segments of $S^1$ transposed into one another
by such a symmetry, then we obtain one and the same $\li_M$- or 
$i^*_M$-route; in this case we should count it only once for any
class of equivalent segments.)

The differential of an $\li_M$-generator is defined exactly in the same way
as in \cite{v90}, \cite{book}, i.e. as the sum of three $\ti_M$-generators
coinciding with it geometrically and taken with appropriate signs
(or only one such generator if this $\li_M$-route has a symmetry of order 3).
\medskip

{\sc An equivalent and more standard description} of 
the term $E_1^{-i,i}$ is as follows
(cf. \cite{BL}, \cite{BN}, \cite{lin}, \cite{kalf} etc.) 

The group $E_1^{-i,i}$ is isomorphic to the space of $G$-valued
functions on the set of (equivalence classes of) $[i]_M$-routes
satisfying the following conditions.

1. {\it Trivial condition}.
For any $[i]_M$-route, having suspicious segments 
such that the corresponding loop in $M$ is contractible,
the value of this function should be equal to 0.

2. {\it Four-term relation.}
Take any $\li_M$-route and realize it by a 
{\it generic} singular knot,
i.e. by an immersion $S^1 \to M$ having 
$i-2$ transverse double points and one triple point,
tangent vectors at which are linearly independent.
This singular knot can be slightly perturbed in 6 different ways,
so that the triple point splits into two double points, and
the $\li$-configuration respected by our knot splits
into a $[i]$-configuration, see fig.~\ref{triple} (= fig. 15 in \cite{v90}).
These 6 perturbations can be divided in a natural 
way into three ordered pairs, numbered by $\ti$-configurations,
coinciding geometrically with the $\li$-configuration;
in fig.~\ref{triple} these pairs are formed by perturbations
1 and 4, 2 and 5, 6 and 3. For any such pair we take the value of our
function on its first member minus the value on the second.
The {\it four-term relation} claims that all 
these three differences of our six terms
should coincide (so that their common value is a
characteristic of the central singular knot; in the previous description
of $E_1^{-i,i}$
this is the value of our function on the corresponding $\li_M$-route).
\medskip

Finally, in the case of nonorientable $M$ (and $G \ne \Z_2$) the group
$E_1^{-i,i}$ consists of $G$-valued functions on the set of
$[i]_M$-routes, satisfying all the same relations
with small modifications  and
additionally taking zero value on all $[i]_M$-routes
such that the corresponding strata of the discriminant do not 
satisfy a certain orientability condition, see \S \ 1.6.3 below.
In the simplest case of $[1]_M$-routes this condition coincides with the
standard coorientability of the corresponding 
stratum of the discriminant variety.
\medskip

{\sc Examples.} 0. Invariants of order $0$ are 
any functions $\pi_0(\OM) \to G$.

1. The unique $[1]$-configuration is the pair of points in the circle.
Any $[1]_M$-route $\alpha \in E_0^{-1,1}$ is 
defined by an immersed circle in $M$ with one
selfintersection point. Its differential $d_0(\alpha) \in \tilde E^{-1,2}_0$ 
is nontrivial if and only if
one of two loops, formed by this circle, defines a zero class
in $\pi_1(M)$.
\medskip

{\sc Proposition 1.} {\it For any oriented 3-dimensional manifold 
$M$, the group of order 1 elements of the group
$H^1(\Omega_f M, \Omega_f M \sm \Sigma)$ is a free Abelian group
whose generators are in the obvious one-to-one correspondence with
the unordered pairs $(\alpha,\beta)\equiv 
(\beta,\alpha) \subset \pi_1(M)$ such that 
$\alpha \ne 0 \ne \beta$, and considered up to simultaneous
conjugations: $(\alpha, \beta) \sim (\alpha', \beta')$
if there is $\lambda \in \pi_1(M)$ such that} 
$\alpha' =\lambda^{-1}\alpha \lambda, \beta' = \lambda^{-1}\beta\lambda$.
\quad $\Box$ 

So, in this theory the non-triviality of the Whitehead link 
\unitlength=1.50mm
\begin{picture}(24,6.5)(-2,2)
\put(15,5){\oval(6,6)[t]}
\put(14.3,5){\oval(4.6,4)[lb]}
\put(16.5,5){\oval(3,4)[rb]}
\put(12.7,4){\oval(4.6,4)[r]}
\put(13.5,4){\oval(5,6)[r]}
\put(13.5,2){\oval(21.4,2)[lb]}
\put(12.7,3){\oval(7.4,2)[lb]}
\put(5.15,2){\oval(4.7,4)[lt]}
\put(6.5,3){\oval(5,2)[rt]}
\put(11.5,5){\oval(11,2)[lt]}
\put(11.7,5){\oval(21.4,4)[lt]}
\put(2.5,5){\oval(3,4)[lb]}
\put(3.65,5){\oval(4.7,4)[rb]}
\put(20.5,4.5){\makebox(0,0)[cc]{$C_1$}}
\put(-1,6){\makebox(0,0)[cc]{$C_2$}}
\end{picture}
in $\R^3$  
can be proved already by an invariant of order 1 (and not of order 3, as in
the usual theory of finite-order invariants). Indeed, 
when we try to deform this link to the trivial one, one of two circles
(say $C_1$) can be considered as unmoved one, hence the
triviality of the link is equivalent to the triviality of the knot
$C_2$ in the manifold $\R^3 \setminus C_1$.
This knot can be obtained from a circle $C_0$, unknotted
and unlinked with $C_1$,
by a deformation, along which it selfintersects only once, and both
loops arising in the instant of selfintersection define nonzero
elements of the group $\pi_1(\R^3 \setminus C_1)$.  
By  Proposition 3 below
the corresponding element of $H^1(\Omega_fM,\Omega_fM \sm \Sigma)$
belongs to the kernel of the obvious map of this group to 
$H^1(\OM)$, and hence defines a knot invariant.
The values of this invariant on knots $C_0$ and $C_2$ differ by $\pm 1$.
\medskip

{\sc Remark.} In the calculation of higher differentials
$d_r$, $r \ge 1$, of the spectral sequence some new generators
of groups $E_*^{-i,i+1}$ can arise, namely, certain 1-dimensional
cohomology classes of spaces of maps $S^1 \to M$ respecting our
$[i]$-configurations, see \S \ 2.4. They can provide some
extra obstructions to the integration process, see \S \ 1.5.

\subsection{Example of a nontrivial finite-order class in $H^1(\OM)$
and some obstructions to the existence of such classes.}

Let $M=S^2 \times S^1$. We construct a 1-parametric family of 
immersed circles in $M$. All these circles consist of two segments, 
the first of which is the same for all circles of the family: it
starts at the north pole of the distinguished sphere
$S^2 \times \{0\}$, finishes at the south pole of the same sphere,
has no self-intersections, and the cyclic coordinate of the
factor $S^1$ grows on it monotonically from $0$ to $4\pi$.
The second segments of these embedded circles will be just
all the meridians in this sphere, joining the south pole back 
to the north one; they (and hence also the entire corresponding
embedded circles in $M$)
are parametrized by the equator circle
in $S^2$. There is exactly one selfintersecting circle 
in this family. 

All curves of our family are piece-wise smooth with only
two breakpoints at the poles; it easy to improve them slightly 
at these points in such a way that they become smooth but 
do not get extra selfintersections. The $[1]_M$-configuration,
arising at the unique instant of self-intersection, defines an element
of $H^1(\OM, \OM \setminus \Sigma)$: indeed, both loops formed by it
are homotopically nontrivial in $M$. The value of this element
on our loop in $\OM$ is equal to $\pm 1$, thus it defines also a nontrivial 
element in $H^1(\OM)$.
\medskip

{\sc Problem.} To construct a similar example inside the trivial
component of $\OM$ (consisting of contractible loops);
may be with a more complicated $M$.
Proposition 3 below shows in particular that it is impossible
if $M$ is closed and $\pi_2(M)=0.$
\medskip

{\sc Proposition 2.} {\it Let $M$ be an arbitrary 3-manifold, and 
$C$ a connected component of 
$\OM$. Suppose that the group $H_1(C,\R)$ contains a basis consisting
of loops $S^1_\lambda \to \OM$ (or, which is the same, of
maps $S^1 \times S_\lambda^1 \to M$) such that for any 
$\lambda \in S^1_\lambda$ the restriction of this map 
on any circle $S^1 \times \{\lambda\}$ is an embedding. 
Then all the elements of $H^1(\OM)$, coming from our spectral
sequence, take zero values on the elements of $H_1(C)$
and hence all elements of groups $E_\infty^{-i,i}$
of this sequence define invariants of knots from this component.}
\quad $\Box$ \medskip

For the ``trivial'' component $C_0$ of $\OM$ (i.e. the component
of the trivial knot) the condition of the previous proposition
follows from the following more standard one.
\medskip

{\sc Proposition 3.} {\it Suppose that  the group $\pi_2(M)$ is generated by 
spheroids $S^2 \to M$ homotopic to embedded spheres. 
Then the condition of the previous proposition for the
component $C_0$ is satisfied, so that
all the elements of $H^1(\OM)$, coming from our spectral
sequence, are trivial in restriction to $H_1(C_0)$.}
\medskip

{\it Proof.} Let us choose a small contractible unknotted parametrized circle 
$\circ \subset M$ for the basepoint in $\Omega_fM$. Given a spheroid
$\varepsilon: S^2 \to M$, denote by $\tau_\varepsilon$ the toroid
$S^1 \times S^1_\lambda \to M,$ which almost everywhere coincides with
the trivial loop $\tau_0$ in $\Omega_fM$ (sending all of $S^1_\lambda$
to the point $\circ$) and only in a small closed disc
$\delta \in S^1 \times S^1_\lambda \setminus (S^1 \times *)$ replacing
$\tau_0$ by such a map that the union of maps $\tau_\varepsilon$
and $\tau_0$ on two copies of $\delta$, glued on their boundaries,
is a spheroid homotopic to
$\varepsilon$.

Any loop in $\Omega_fM$, i.e. a continuous map $S^1 \times S^1_\lambda \to M$,
is homotopic to a sequence of loops, the first of which is the family
of small unknotted embedded circles, moving along a path in $M$, and
all other are toroids $\tau_{\varepsilon_j}$, where $\varepsilon_j$ are some
basis elements of $\pi_2(M)$. Thus we need only to check the
conditions of Proposition 2 for any toroid $\tau = \tau_\varepsilon$,
where the spheroid $\varepsilon$ is homotopic to an embedding.

By definition, there is 
a small neighborhood
$I$ of the distinguished point $* \in S^1_\lambda$ such that for $\lambda \in I$
the corresponding loops $\tau(S^1 \times \{\lambda\})$ coincide
with $\circ$. 
Denote by $D = D_N \cup D_S$ the union of two discs in $S^2$
bounded by small polar circles, and consider the spheroid $s: S^2 \to M,$
coinciding in $S^2 \setminus D$ with $\tau$ 
after the standard identification
\begin{equation}
\label{polar}
S^1 \times (S^1_\lambda \setminus I) \sim S^2 \setminus D,
\end{equation}
and in $D$ coinciding with two small embedded discs contracting 
the circle $\circ$ in $M$.

Let $\theta: [0,1] \times S^2 \to M$ be a homotopy of this spheroid,
moving it to an embedding of $S^2$. Since for a generic map $S^2 \to M$
the set of singular points is finite, we can assume that during
this homotopy the restrictions of all maps $\theta(t,\cdot)$
on the polar zone $D$ are immersions
(and thus also embeddings since this zone can be chosen
arbitrarily small). Then consider the
toroid $S^1 \times S^1_\lambda \to M$, whose restriction on the
cylinder $S^1 \times (S_\lambda^1 \setminus I)$ 
coincides via the identification (\ref{polar})
with the spheroid $\theta(1, \cdot)$, and the image of the cylinder
$S^1 \times I$ coincides with the union of 
traces $\theta([0,1] \times \partial D)$
of polar circles during the homotopy (glued on their
common circle $\theta(0 \times \partial D_N) \equiv 
\theta(0 \times \partial D_S)$). 
This toroid satisfies conditions of Proposition 2.
\quad $\Box$ \medskip

{\sc Corollary.} There exists a non-trivial theory of finite-order invariants
in reducible 3-manifolds, cf.~\cite{kalf}. For instance, $S^2 \times S^1$
is reducible but satisfies conditions of Proposition 3.
\medskip

I am sure that experts in the low-dimensional
topology can prove much stronger statements also proving
the triviality of $H^1(\OM)$-classes of this kind, cf.~\cite{kalf}.

\subsection{Integration of elements of $E_1^{-i,i}$ to knot invariants.}

In this subsection we assume that $M$ is a three-dimensional
oriented manifold. Given an element $\gamma$ of the corresponding group
$E_1^{-i,i}$, can it be continued to a knot invariant of order $i$?
Similarly to \cite{v90}, \cite{BL}, this question can be reduced to a
sequence of systems of linear equations, whose unknowns correspond to 
topological types of singular knots in $M$ with $<i$ selfintersections,
cf. \cite{lin}, \cite{kalf}. In this subsection we reformulate this
condition back in terms of spectral sequences: $\gamma$ should belong
to the subgroup $E_{i+1}^{-i,i} \subset E_1^{-i,i}$, see \cite{v90}.
This allows us to write explicitly the group of additional obstructions
to this continuation, arising from the topology of the manifold $M$,
see formula (\ref{circ}) below.

\subsubsection{Actuality table.}

Similarly to \cite{v90}, any  
knot invariant 
of a finite order $i$ can be encoded by
an {\it actuality table}, having $i+1$ levels $0,1,2, \ldots,i$.
The $j$-th level consists of cells, corresponding to all
$[j]_M$-routes.

In any such cell we draw a generic immersion $S^1 \to M$ representing
this $[j]_M$-route (i.e. having exactly $j$ transverse selfintersection points
and no other singularities). This immersion is an accessoir of the
table itself and does not depend on the invariant. To describe the
invariant, we write in the cell a number (or, more generally, the element
of the coefficient group $G$), namely, the index
of $j$-th order of the corresponding immersion, see \S \ 1.2.  

The calculation of the value of an invariant consists in the same inductive
process as in \cite{v90}, \cite{book}. Namely, 
we join our knot by a generic path
in $\Omega_fM$ with the {\it distinguished} knot from the same
homotopy class in $\Omega_fM$ (i.e., to the knot drawn in the corresponding
cell of the $0$-th level of the actuality table).
Such a path has only finitely many intersection points with the 
variety $\Sigma$ at its points, corresponding to immersions 
$S^1 \to M$ with one transverse self-intersection.

The value of the invariant at our knot is equal to its value
at the distinguished one plus the sum of 
1-st order indices of these immersions,
taken with coefficients $\pm 1$,
equal to signs of corresponding local
surgeries, see the last definition in \S \ 1.2.
To calculate these indices, we
join these immersions to distinguished ones (given in the first level 
of the table) by arbitrary generic paths in $\Sigma$ and count all the
points of the set of transversal self-intersection of $\Sigma$, which we meet
along these paths, etc. This process stops on the level $i$, because
by the characteristic property (see \S \ 1.2) any invariant of order $i$
defines the same indices of order $i$ for all generic immersions from the
same $[i]_M$-route. (In particular, we do not need to draw pictures
in the actuality table at the top level.)

If we do not fix values of the invariant at the $0$-th level
of the table, then we get not an invariant, but just an element
of the group $H^1(\Omega_fM, \Omega_fM \sm \Sigma)$.  

Our spectral sequence (more precisely, its restriction on 
cells $E_r^{-j,j}, E_r^{-j,j+1}$, responsible for the
$0$-dimensional cohomology of the space of knots) is a method 
of calculating all actuality
tables such that this algorithm works and does not lead to
a contradiction.

It starts with any element $\gamma \in E_1^{-i,i}$, see \S \ 1.3. 
This element defines an upper ($i$-th) level of the actuality table:
into any cell, corresponding to an $[i]_M$-route, we put the value of
the function $\gamma$ at this route. Then, exactly as in 
\cite{v90}, we fill in the table from top to bottom. We can fill in
all levels $i-1, i-2, \ldots, i-r+1$ if and only if $\gamma$
belongs to the subgroup $E_r^{-i,i}$ of $E_1^{-i,i}$. Unlike
\cite{v90}, this subgroup can be proper, see \S \ 1.5.4 below.

\subsubsection{The short spectral sequence.}

Exactly as in \cite{v90}, we factorize our general spectral sequence
(described in \S \ 2) through some elements which surely do not 
contribute to the calculation of $0$-di\-men\-si\-o\-nal cohomology of the space
of knots in $M$, and obtain the {\it short spectral sequence}
$E_r^{p,q}$, $r \ge 0$, whose non-trivial groups lie only on 
two half-lines with $p<0$ and $p+q = 0 $ or $p+q=1$, and coincide with
the corresponding groups of the main spectral sequence on the first
of these lines. In the rest of the present section we deal only with
this short spectral sequence.

The groups $E_0^{-i,i}$ are already described in \S \ 1.3, let us
describe $E_0^{-i,i+1}$.

For any $[i]_M$-route $I$ denote by $\{I\}$ the space of its realizations,
i.e., of pairs of the form \{an $[i]$-configuration of the corresponding
equivalence class; a map $S^1 \to M$ respecting this configuration\}.

The group $E_0^{-i,i+1}$ is the direct sum of the group $\tilde E_0^{-i,i+1}$,
described in \S \ 1.3, and the group
\begin{equation}
\label{circ}
\breve E_0^{-i,i+1} \equiv \prod_I H^1(\{I\},G),
\end{equation}
multiplication over all $[i]_M$-routes $I$.

\subsubsection{The differential $d^1: E_1^{-i,i} \to E_1^{-i+1,i}$.}

Suppose that we have an element $\gamma$ of $E_1^{-i,i}$, i.e.
a $G$-valued function  on the set of $[i]_M$-routes,
satisfying the basic relations described in \S \ 1.3. 
To extend it to a knot invariant, we need to calculate all
its higher differentials $d^1(\gamma) \in E_1^{-i+1,i}$,
$d^2(\gamma) \in E_2^{-i+2,i-1}, \ldots,$ 
$d^{i-1}(\gamma) \in E_{i-1}^{-1,2}$, $d^i(\gamma) \in 
E_i^{0,1}$, and to prove that all these differentials are 
trivial; if we prove all these conditions but the last one,
then we extend our element $\gamma$ only to an element of the
group $H^1(\Omega_f M, \Omega_f M \sm \Sigma)$.

In this subsubsection we describe explicitly the first of these conditions.
\medskip

The map $d^1$ splits into the sum of two operators
$\breve d^1: E_1^{-i,i} \to \breve E_1^{-i+1,i}$
and $\tilde d^1: E_1^{-i,i} \to \tilde E_1^{-i+1,i},$
where 
$\breve E_1^{-i+1,i} \equiv \breve E_0^{-i+1,i}$,
$\tilde E_1^{-i+1,i} \equiv \tilde E_0^{-i+1,i}/d^0(E_0^{-i+1,i-1})$,
and we need to check both conditions
$\breve d^1(\gamma) =0$, $\tilde d^1(\gamma) = 0$.
Let us describe these conditions.

Consider any $[i-1]_M$-route $I$, and any 1-cycle
$l$ in the manifold $\{I\}$. We can realize it by a smooth generic path,
only finitely many times intersecting the set of immersions,
having $i$ transverse self-intersections, i.e. defining
points of $[i]_M$-routes. At any such point we take the 
value of our function $\gamma$ on this $[i]_M$-route, and multiply
it by the coefficient $1$ or $-1$ depending on the direction in which
we traverse this stratum (i.e. on the sign of the corresponding
local surgery close to the additional $i$-th crossing point, see \S \ 1.2).
The sum of all such values taken with these coefficients depends only
on the class of the path $l$ in the group
$H_1(\{I\},G)$. We denote this sum by
$\langle \breve d^1(\gamma), l\rangle$ and thus 
define the element 
$\breve d^1(\gamma) \in H^1(\bigcup_I \{I\};G).$

If this element is non-trivial, then $\gamma$ is not equal to the 
$i$-th index of any invariant of order $i$.

Now suppose that $\breve d^1(\gamma)=0$. Then we can define a locally
constant function $A_{i-1}$ on the regular part 
of the manifold $\{I\}$ (i.e. on 
its part corresponding to immersions with exactly $i-1$ selfintersections)
in such a way that the difference of its values 
on two sides of any hypersurface, corresponding to any $[i]_M$-route,
coincides with the value of $\gamma$ at this route. Such a 
locally constant function $A_{i-1}$ is 
defined by $\gamma$ almost uniquely, only up to addition of functions, 
which are constant on any manifold $\{I\}$ (i.e., up to
elements of the group $E_0^{-i+1,i-1}$). 

The second condition
$\tilde d^1(\gamma)= 0$, which this function $A_{i-1}$
should satisfy, consists of following two subconditions: 

a) given any generic immersion 
with $i-3$ transverse double crossings and one triple point, the values of
$A_{i-1}$ at all 6 local moves 
of the triple point,
decomposing it into a pair of double points, satisfy
the 4-term relations, see fig.~\ref{triple};

b) given any generic map $S^1 \to M$ with $i-2$ transverse self-intersections
and one cusp point, the value of $A_{i-1}$ at its local move, replacing the
cusp by a self-intersection point, is equal to $0$.

It is sufficient to check this second condition close to one generic point
of any $\langle i-1 \rangle_M$-route or $(i-1)^*$-route:
if it is satisfied for some choice of such points (and the first 
homological condition $\breve d^1(\gamma)=0$ also holds),
then it will be satisfied automatically
at all other points of the same $\langle i-1 \rangle_M$-route
or $(i-1)^*$-route.

This condition can be easily identified as the 
triviality of a certain element $d^1(\gamma)$ of the 
quotient group $\tilde E_1^{-i+1,i} \equiv
E_0^{-i+1,i} / d^0(E_0^{-i+1,i-1})$.

Suppose that  
$d^1(\gamma)=0$, i.e. there exists a locally constant function
$A_{i-1}$ satisfying all above conditions. Then we say that $\gamma$ belongs
to the group $E_2^{-i,i}$, and are able to fill in the $(i-1)$-th level of
the actuality table: in any cell, containing a generic point of
the $[i-1]_M$-route, we put the number (or element of $G$) equal
to the value of the function $A_{i-1}$ at this point. 

The choice of this function $A_{i-1}$ is not unique (if exists): 
it is defined by $\gamma$ up to 
addition of arbitrary
elements of the group $E_1^{-i+1,i-1}$.

\subsubsection{An example of non-degenerating spectral sequence.}

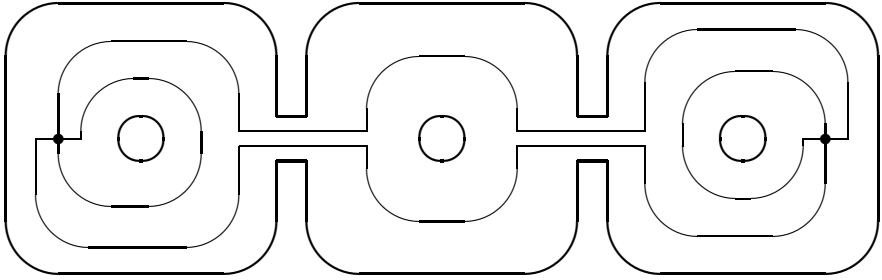
\begin{figure}
\unitlength=1.00mm
\special{em:linewidth 0.4pt}
\linethickness{0.4pt}
\begin{picture}(130,40)
\thicklines
\put(18,15){\oval(36,30)[b]}
\put(18,21){\oval(36,30)[t]}
\put(58,15){\oval(36,30)[b]}
\put(58,21){\oval(36,30)[t]}
\put(98,15){\oval(36,30)[b]}
\put(98,21){\oval(36,30)[t]}
\put(18,18){\oval(6,6)}
\put(58,18){\oval(6,6)}
\put(98,18){\oval(6,6)}
\put(0,15){\line(0,1){6}}
\put(36,15){\line(1,0){4}}
\put(36,21){\line(1,0){4}}
\put(76,15){\line(1,0){4}}
\put(76,21){\line(1,0){4}}
\put(116,15){\line(0,1){6}}
\thinlines
\put(31,17){\line(1,0){17}}
\put(31,19){\line(1,0){17}}
\put(68,17){\line(1,0){17}}
\put(68,19){\line(1,0){17}}
\put(4,18){\line(1,0){6}}
\put(106,18){\line(1,0){6}}
\put(112,19){\line(0,-1){1}}
\put(4,17){\line(0,1){1}}
\put(18,18){\oval(16,16)[t]}
\put(19,19){\oval(24,24)[t]}
\put(16.5,19){\oval(19,20)[b]}
\put(17.5,17){\oval(27,27)[b]}
\put(58,17){\oval(20,20)[b]}
\put(58,19){\oval(20,20)[t]}
\put(98.5,19){\oval(27,27)[t]}
\put(97,17){\oval(24,24)[b]}
\put(98,18){\oval(16,16)[b]}
\put(99.5,17){\oval(19,20)[t]}
\put(7,18){\circle*{1.33}}
\put(109,18){\circle*{1.33}}
\end{picture}
\caption{A non-integrable singular knot.}
\label{int}
\end{figure}

Let $M$ be the connected sum of three copies of $S^2 \times S^1$.
A planar outline of $M$ is shown in fig.~\ref{int} by the domain
with three holes, bounded by thick curves.
Consider the loop $S^1 \to M$ with two self-intersections, 
shown in fig.~\ref{int} by the thin line. Its chord diagram 
is {\it trivial}, i.e., consists of two non-crossing chords.

\medskip
{\sc Definition.} Two $[2]_M$-routes with non-crossing chords
are {\it neighbors}, if there exists an immersion $S^1 \to M$ with
unique generic triple point  such that some two of three
its perturbations, shown in fig.~\ref{triple} and respecting 
trivial chord diagrams,
belong to these $[2]_M$-routes.
Two $[2]_M$-routes are {\it related}, if there exists a chain
of $[2]_M$-routes, joining them, any two neighboring members of which
are neighbors. A $[2]_M$-route is {\it marginal} if
one of two its suspicious loops (see \S \ 1.3)
is contractible (so that any element of the group $E_1^{-2,2}$
should take zero value on the corresponding $[2]_M$-route).
\medskip

{\sc Lemma.} {\it Among relatives of the $[2]_M$-route, represented by the 
curve from fig.~\ref{int},
there are no marginals.} 
\medskip

Indeed, the subgroup in $H_1(M)$, generated by cycles, lying in
an immersed circle, is the same for all its relatives. For the initial
curve from fig. \ref{int} this subgroup is of rank 3,
and for any marginal $[2]_M$-route of rank at most $2$. \quad $\Box$
\medskip

Now define the $\Z$-valued function $\gamma$ on the set of all $[2]_M$-routes,
which takes value 1 at all relatives of the route represented by the
curve from fig. \ref{int}, and zero value at all other $[2]_M$-routes.
This function satisfies both conditions from \S \ 1.3, and hence belongs
to the subgroup $E_1^{-2,2} \subset E_0^{-2,2}$.

However $\breve d^1(\gamma) \ne 0$. Indeed, consider 
the horizontal segment in the picture of our curve in the rightmost
copy of $S^2 \times S^1$, joining two points of the sphere 
$S^2 \times \{0\}$. We can suppose that these points are 
poles of this sphere, and the segment is its distinguished meridian.
Consider the family 
of curves $C_\alpha$, $\alpha \in [0,2\pi]$, 
coinciding with the one from 
fig. \ref{int} everywhere outside this segment 
and replacing it by all other meridians; the parameter 
$\alpha$ of this family is 
the cyclic coordinate $\alpha$ of the equator in $S^2$.

This family defines a 1-cycle in the corresponding   
$[1]_M$-route, and its intersection index with the 
cycle $\gamma$ obviously is equal to $\pm 1$, in particular
$\gamma$ defines a non-zero class 
in the group $H^1$ of this $[1]_M$-route. 
Thus $\breve d^1(\gamma) \ne 0.$

\subsubsection{Differentials $d^2, d^3$ etc.}

Now suppose that 
we already have calculated (and defined) all differentials 
$d^s:E_s^{-j,j} \to E_s^{-j+s,j-s}$, $s < r$, for any $j \le i$;
let $\gamma$ be an element of the subgroup
$E_r^{-i,i} \subset E_1^{-i,i}$. This means in particular that we  
can fill in all levels $i-1, i-2, \ldots,
i-r+1$ of the actuality table with the upper level $\gamma$. The
indices of this table are then determined by $\gamma$
up to addition of similar tables corresponding to all elements
of groups $E_{r-1}^{-i+1,i-1}, E_{r-2}^{-i+2,i-2}, \ldots,
E_1^{-i+r-1,i-r+1}$.

Moreover, in this case we can define similar indices 
(of corresponding orders)
for {\em all} immersions $S^1 \to M$, having exactly
$i, i-1, i-2, \ldots, i-r+1$ transverse self-intersection points.
Then for any manifold $\{I\}$, where $I$ is a $[j]_M$-route with
$j \ge i-r+1$,
these indices form a locally constant function on its open
submanifold, consisting of immersions with exactly $j$ transverse
self-intersections. Let us fix such a collection 
of locally constant functions $A_j$
and the actuality table, representing them, in which only levels
$i, i-1, \ldots, i-r+1$ are filled in.
We call this partially completed table the {\it tentative actuality table}.

This collection of locally constant functions (or, 
equivalently, the representing 
them tentative actuality table) is called $1$-{\it integrable,} if there 
exists a locally constant function on the union of regular subsets
of all manifolds $\{I\}$ for all $[i-r]_M$-routes $I$, such that
a) the difference of its values at two sides of any hypersurface
in $\{I\}$, consisting of generic immersions with $i-r+1$ self-intersections,
is equal to the value at the corresponding piece of this
hypersurface of its index $A_{i-r+1}$, encoded in the 
existing part of table, and
b) close to any generic immersion, respecting an 
$\langle i-r \rangle$-configuration (respectively, $(i-r)^*$-configuration), 
all corresponding 4-term relations (respectively, the trivial
relation) are satisfied.

Exactly as in \S \ 1.5.3, the obstruction to the existence of such a function
is just an element of the group $E_1^{-i+r,i-r+1}$.
Denote by $D^r(E_r^{-i,i})$ the subgroup of this group, generated by
all such obstructions over all elements $\gamma \in E_r^{-i,i}$
and all their admissible extensions to collections of locally constant
functions $A_j$, $j=i-1, \ldots, i-r+1$. 

The condition 
$\gamma \in E_{r+1}^{-i,i}$
means that this obstruction is trivial for at least one 
choice  of the tentative actuality table with 
upper level equal to $\gamma$. To check this condition, we need
to calculate this obstruction for an arbitrary such tentative table,
and check whether it belongs to the subgroup generated by
similar subgroups

\begin{equation}
\label{relations}
d^1(E_1^{-i+r-1,i-r+1}), D^2(E_2^{-i+r-2,i-r+2}), \ldots,
D^{r-1}(E_{r-1}^{-i+1,i-1}).
\end{equation}

\noindent
Here is one more way to say the same: we denote by $E_r^{-i+r,i-r+1}$
the quotient group of $E_1^{-i+r,i-r+1}$ by the subgroup generated
by all subgroups (\ref{relations}), denote by $d^r(\gamma)$ the class
of our obstruction in this quotient group, and check
whether it is trivial or not.

If yes, then we a) change the levels
$i-1, \ldots, i-r+1$ of our tentative actuality table by the elements of
an arbitrary $1$-integrable one with the same leading 
term $\gamma$, and b) fill in
the $(i-r)$-th level of this new table, putting in any cell, corresponding to
a $[i-r]_M$-route, the value, which any locally constant 
function $A_{i-r}$ on this route,
satisfying the above integration conditions (and defined by this 
$1$-integrable tentative table), takes at the
immersion depicted in this cell.

Thus the inductive step of the calculation
(and definition) of our spectral sequence is completed.
\medskip

{\sc Remark.} We could replace the group (\ref{circ}) by certain its
subgroup. Namely, for an $[i]_M$-route $I$ denote by $h(I,G)$
the subgroup in $H_1(\{I\},G)$ spanned by such loops in $\{I\}$,
all whose points are generic immersions with exactly $i$ selfintersections.
Then the group 
$\breve E_0^{-i,i+1}$ could be redefined as
$\prod_I(\mbox{Ann}\ h(I,G)) \subset \prod_I H^1(\{I\},G)$. Indeed, all our
additional obstructions to integrability take zero values on elements of 
$h(I,G)$.

\subsection{The case of non-orientable manifolds.}

If the three-dimensional manifold $M$ is non-orientable, 
then many natural strata of (the
resolution of) the discriminant turn out to be non-(co)orientable,
and hence can participate only in the construction of (mod 2)-invariants.
However, if $M$ has a sufficiently complicated fundamental group, then
many strata are still orientable and define integer invariants and
homology classes.

\subsubsection{First example.}

Consider the main stratum of the 
discriminant $\Sigma \subset \Omega_f M$, corresponding
to some $[1]_M$-route (i.e. an 
irreducible component of the set of maps $S^1 \to M$ 
with a transverse self-intersection).
Let $\phi$ be a generic representative of this stratum, and 
$\alpha$ and $\beta$ the classes of two corresponding
loops in the group $\pi_1(M)$ with basepoint at the self-intersection point.
\medskip

{\sc Proposition 4.} {\it This stratum of the discriminant is non-orientable
(i.e. its intersection with any sufficiently large approximating
space $\Gamma_U$ is)
if and only if there exists an element $\lambda \in \pi_1(M)$
such that 

a) $\lambda$ destroys the orientation of $M$: 
$\langle w_1(M), \lambda \rangle \ne 0$;

b) the conjugation operator $T_\lambda:\pi_1(M) \to \pi_1(M)$,
$T_\lambda(\cdot) = \lambda^{-1}(\cdot)\lambda$, either preserves
both elements $\alpha$ and $\beta$, or permutes them.}
\quad $\Box$ \medskip

Consider the connected sum $M^2 = K \# K$ of two copies 
of the Klein bottle, and an immersed curve 
in it with unique self-intersection on the ``neck'' of the connected sum,
such that any of two obtained loops lies in its own summand of
the connected sum and defines in it a basic loop destroying its orientation.
Using the obvious identification $M^2 = M^2 \times \{0\}$, we will
consider this curve as a loop in the manifold $M^3 = M^2 \times S^1$.
It follows from the van Kampen theorem, that an element
$\lambda \in \pi_1(M^3),$ satisfying the conditions of the 
previous proposition,
does not exist. Hence our stratum defines a class in the integer cohomology
group $H^1(\OM, \OM \setminus \Sigma)$. 
\medskip

It is easy to prove that the condition of Proposition 2 
is satisfied for the containing this curve component of
$\OM$, hence this class defines a knot invariant.
There are infinitely many knots distinguished
by this invariant. Namely, 
let $\nu$ be the ``neck'' cylinder connecting two summands $K$ of $M^2$.
Everywhere outside $\nu \times S^1$ our knots coincide with our
immersed circle in $M^2 \times \{0\}$, and in $\nu \times S^1$
they coincide with 2-string braids going from one boundary component to
the other and twisted arbitrarily many times.

\subsubsection{Coding order $i$ invariants in non-orientable
3-manifolds and computation of their values.}

Consider any immersion $\phi: S^1 \to M^3$ with $i$ 
distinguished transverse double crossings. As usual, we can resolve all these
crossings in $2^i$ locally different ways
so that they become nonsingular knots.
These $2^i$ resolutions can be obviously partitioned into two groups in
such a way that any two neighboring resolutions (i.e., two resolutions, 
obtained from one another by one local surgery of fig.~\ref{cp}) belong to
different groups. 
\medskip

{\sc Definition.} The {\it supercoorientation} (or simply 
$s$-{\it orientation}) of the containing $\phi$ $[i]_M$-route $I$ 
is the simultaneous choice of one of 
these two groups close to all points of the manifold $\{I\}$,
depending continuously of these points.
\medskip

For singular knots in oriented manifolds these $s$-orientations are
defined canonically, see \S \ 1.2.
The problem of deciding whether an $[i]_M$-route in a non-orientable manifold
satisfies this condition or not, is an independent problem,
related in particular with its symmetry properties. 
If $\pi_1(M)$ is sufficiently complicated and 
the loops in $M$, lying in the representing this $[i]_M$-route curve,
are ``sufficiently independent''
in (the set of conjugacy classes of) $\pi_1(M)$, then this route 
should be $s$-orientable.
\medskip

{\sc Definition.}
The {\it global stratum}
(or simply {\it stratum}) of
$\Sigma$ corresponding to an $[j]_M$-route is an irreducible
component of the set of maps $S^1 \to M$, respecting this route 
and having transverse crossings at all its $j$ double points.
(This transversality condition is not very restrictive, 
indeed, the set of maps not satisfying it
has codimension 2 in the set of all maps respecting this route.)
Such a stratum contains an open smooth subset, consisting of maps, 
having no extra
singular or multiple points; its path-components are called {\it small
strata} or {\it pieces} of the global stratum.
\medskip

Any knot invariant can be 
extended in the almost standard way (cf. \S \ 1.2) to a function 
on immersions with arbitrarily
many double transverse crossings and no other singularities. More precisely, 
this extension takes values
on pairs of the form \{such an immersion, 
a choice of the local $s$-orientation of its route\}.
It is equal to the sum of values of our invariant over all resolutions
from the chosen group minus the similar sum over the remaining group;
in particular
it changes the sign if we change the $s$-orientation.
The invariant is of order $i$ if and only if its extension to all 
singular immersions with $>i$ crossings is equal to 0.

Any invariant of order $i$ is encoded by the actuality table,
having $i+1$ levels $0,1, \ldots, i$. The cells of this table
on level $j \le i$ are numbered by 
$[j]_M$-routes.

In any such cell we draw an immersion $\phi$, representing 
this $[j]_M$-route $I$, having exactly $j$
transverse crossings, and supplied with a certain
$s$-orientation of the corresponding small stratum of $I$ 
(these data do not 
depend on the invariant). To specify the invariant, we 
write in this cell  a number: the
index of order $j$ of our invariant at this ($s$-oriented) 
singular immersion, cf. \cite{v90}, \cite{book}.
\medskip

{\sc Remark.} We need to fill in the cells even for
not $s$-orientable $[j]_M$-routes, however if already the corresponding 
small stratum is not $s$-coorientable, then the corresponding
index will be equal to zero.
\medskip

The calculation
of values of invariants
on knots consists essentially in the same inductive process
as in \S \ 4.3 of \cite{v90} (see also \S 1.5.1 above),
only with following modifications.

Suppose that we go along a smooth path in our stratum, 
consisting of immersions with
$j$ crossings, and at some instant traverse the stratum of
$(j+1)$-crossed curves (i.e., at that instant our curve has the
$(j+1)$-th selfintersection point).
We need to compare three $s$-orientations:
these of our $[j]_M$-route at some its point $a_+$
before traversing, after it (at the point $a_-$ of the same global
stratum), and the
$s$-orientation of the $[j+1]_M$-route at the very point $a_0$ of traversing.
Set $\alpha=1$ if two first $s$-orientations are compatible in obvious way
(i.e. the same perturbations of first $j$ crossings will belong
to the chosen groups independently on what is happening close to
the $(j+1)$-st one). Otherwise set $\alpha=-1$.

Set $\beta=1$ if the perturbations from the chosen group close
to the $[j]$-stratum at the point $a_+$ belong also to the group chosen
in correspondence with the $s$-orientation of the $[j+1]_M$-route
at the point $a_0$. Otherwise set $\beta=-1$.

Finally, the value of the index (of order $j$) of our invariant 
at the point $a_+$ of the $[j]_M$-route
is equal to that at the point $a_-$, taken with the coefficient $\alpha$,
plus the value of the index of order $j+1$ at the point $a_0$ of 
the $[j+1]_M$-route, taken with the coefficient $\beta$.

\subsubsection{Basic relations for non-orientable manifolds.}

The basic relations defining the group $E_1^{-i,i}$ of the spectral
sequence (see \S \ 1.3) should be slightly modified in the case of 
non-orientable $M$. 

First, our $G$-valued function on the set of $[i]_M$-routes
should vanish on all  non-$s$-orientable routes. 
The trivial relation stays unchanged, and in the 4-term
relation our 6 perturbations should be taken with certain signs,
depending on their $s$-orientations. To define them we need the following
notion of the $s$-orientability of $\li_M$-routes.

Consider a {\it generic} map $\phi$ realizing some
$\li_M$-route. In its small neighborhood the pair 
$(\Gamma_U,\Sigma \cap \Gamma_U)$
is diffeomorphic to the direct product of the $(\dim \Gamma -i-1)$-dimensional
real space and the pair ($\R^{i+1}$, the union of coordinate 
hyperplanes in $\R^{i+1}$). In
particular in this neighborhood the discriminant locally 
separates the space of knots into $2^{i+1}$ octants.

Let us divide all these octants into two groups in
such a way that any two neighboring octants belong to different groups.
\medskip

{\sc Definition.} The (local)
$s$-orientation of 
the $\li_M$-route, containing $\phi$, 
is a choice of one of these two groups.

The {\it global} $s$-orientation of the $\li_M$-route
is its simultaneous local $s$-orientation at all its 
generic points, depending continuously of these points and compatible
in the obvious way close to  generic points, respecting 
$\langle i+1 \rangle$-configurations. (I.e., if we move along a path in the
$\li_M$-route and traverse the set of maps, having 
an additional selfintersection point, not participating in the
definition of the $\li_M$-route, then local $s$-orientations,
defined in the terms of resolutions of multiple points, participating in
this definition, should not remark this traversing.)

The $s$-orientation (local or global) of an 
$ \langle \tilde i \rangle_M$-route is just
the $s$-orientation of the corresponding $\li_M$-route.

The  $s$-orientation of a $i^*_M$-route $\{I\}$ at its point, having $i-1$
transverse self-in\-ter\-sec\-ti\-ons, is any  
local $s$-orientation of the $[i-1]_M$-route, 
obtained from $I$ by forgetting about its singular point.
\medskip

Let us fix any local $s$-orientation of our $\li_M$-route at
its generic point $\phi$. The union of all $[i]_M$-routes is
represented close to $\phi$ by 
6 locally different components, see fig.~\ref{triple}. 
Supply any of these components 
with a sign, equal to $1$ or $-1$
depending on whether the restriction of this $s$-orientation of the
$\li_M$-route at $\phi$
onto the set of $2^i$ local resolutions of singular knots from this
component coincides with the own $s$-orientation of the corresponding
$[i]_M$-route or not.

Then three {\bf sums} of values of our $G$-valued function
on perturbations 1 and 4 (respectively, 2 and 5, respectively, 3 and 6, see
fig.~\ref{triple}),
taken with these coefficients, should coincide.

\subsubsection{Higher obstructions to the integration.}

Given an element $\gamma \in E_1^{-i,i}$ (i.e., a function on the
set of $s$-oriented $[i]_M$-routes, satisfying the relations
from \S \ 1.6.3), its integration to an order $i$ knot invariant with
upper level $\gamma$ (in particular obstructions to the existence
of such an invariant) can be formulated in terms of the short spectral sequence,
generalizing that from \S \ 1.5.

Again, its non-trivial groups $E_r^{p,q}$ lie on only two lines
$p+q=0$ and $p+q=1$, i.e. are of the form $E_r^{-i,i}$
or $E_r^{-i,i+1}$, $i \ge 0$.

Its group $E_0^{-i,i}$ is the space of $G$-valued functions on
the space of $s$-oriented $[i]_M$- and $\li_M$-routes, taking 
opposite values on any route supplied with opposite $s$-orientations.

The group $E_0^{-i,i+1}$ is the sum of two groups
$\tilde E_0^{-i,i+1}$ and $\breve E_0^{-i,i+1}$.
The first of them is generated by all $s$-oriented 
$\langle \tilde i \rangle_M$- and $i^*$-routes.
The second is defined by
\begin{equation}
\label{nonor}
\breve E_0^{-i,i+1} \equiv \prod_I H^1(\{I\}, sG),
\end{equation}
summation over {\em all} $[i]_M$-routes $I$, where $sG$ is the local
system of groups, locally isomorphic to $G$ and such that the monodromy
over a loop in $\{I\}$, destroying the $s$-orientation, 
acts in the fibre as multiplication by $-1$.

The operator $d^0$ acts from $E_0^{-i,i}$ to the first summand
$\tilde E_0^{-i,i+1}$, and its kernel is (naturally isomorphic to)
the group $E_1^{-i,i}$ described in \S \ 1.6.3. The construction of
forthcoming operators $d^r$ essentially repeats that from
\S \ 1.5.

\medskip
{\sc Remark.} All considerations and events from sections 1.6.2---1.6.4 
are valid if $M$ is orientable and coincide then with
their standard versions, see \cite{v90}, \cite{BN}.

\subsection{Functoriality of spectral sequences.} 

If $M'$  is a submanifold of $M$ (of the same dimension), then
there appears the natural homomorphism of our spectral sequences,
$E_r^{p,q}(M) \to E_r^{p,q}(M')$. Indeed, the space of maps
$S^1 \to M'$ is an open subset in the space of maps $S^1 \to M$.
(In the framework of finitedimensional approximations from 
\S \ 1.1, for the approximating 
set of the space $\OM'$ we can take the subset in $\Gamma_U,$
consisting of maps, whose images belong to $\tau^{-1}(M')$.)
This embedding induces the restriction homomorphism
from the Borel--Moore homology group of 
the discriminant set of the former space to that for the latter
one. This homomorphism can be extended naturally 
to spaces of resolutions of these discriminants
and to any terms of their natural filtrations, thus inducing 
a homomorphism of spectral sequences, see \S \ 2.3. The explicit form of
these spectral sequences 
implies the following theorem.
\medskip

For any $[i]_M$-route $I$ denote by $\{I_{M'}\}$
the subspace in $\{I\},$ formed by maps, whose images
belong to $M'$.
\medskip

{\sc Theorem 1.} {\it Let $M' \subset M$ be two three-dimensional 
oriented manifolds, and suppose that the identical embedding 
$M' \to M$ induces 

a) an isomorphism $\pi_1(M') \to \pi_1(M)$, and

b) for any $[i]_M$-route $I$ an epimorphism 
(respectively, isomorphism) $H_1(\{I_{M'}\}) \to H_1(\{I\})$.

Then the group of finite-order invariants of knots  
(or, more generally, $d$-component links with any fixed $d$) in the manifold
$M'$ is canonically isomorphic to a quotient group of the similar group
of invariants of knots or links in $M$ (respectively, to all 
this group).}
\medskip

Indeed, our condition a) implies that for any $i$ our embedding induces
the natural isomorphism $E_1^{-i,i}(M) \to E_1^{-i,i}(M')$,
and condition b) implies that all the forthcoming maps
$E^{-i,i}_r(M) \to E^{-i,i}_r(M')$, $r>1$, are epimorphic (respectively,
isomorphic), so that the limit homology map
$H^1(\OM,\OM \setminus \Sigma) \to 
H^1(\OM', \OM' \setminus \Sigma)$ also is epimorphic (respectively, 
isomorphic). Moreover, condition
b), applied to the $[0]_M$-routes (i.e. path-components of spaces
$\OM$ and $\OM'$) implies that the kernel of the map
$H^1(\OM', \OM' \setminus \Sigma) \to H^1(\OM')$ is a quotient group of the 
similar kernel for $M$. \quad $\Box$ \medskip

The isomorphism theorem from \cite{lin} follows immediately
from this one, see \cite{congr}. Indeed, we can take $M' = \R^3$, then 
our conditions a) and b) will be satisfied for all 2-connected 3-manifolds
$M$. \medskip

Here is a slightly more general statement.
\medskip

{\sc Theorem 1$'$.} {\it Let $M' \subset M$ be two three-dimensional
oriented manifolds, such that

a) for any $[i]_M$-route $I$ the space $\{I_{M'}\}$ consists of
at most one path-component (respectively, of exactly one), and

b) for any $I$ such that $\{I_{M'}\}$ is non-empty,
the map $H_1(\{I_{M'}\}) \to H_1(\{I\})$, 
induced by the identical 
embedding, is epimorphic (respectively, isomorphic). 

Then the graded group of finite-order invariants of knots in  
$M'$ is naturally isomorphic to a quotient group of the similar 
group of invariants of knots in
$M$ (respectively, to entire this group).}  
\medskip

Indeed, condition a) ensures that the map 
$E_1^{-i,i}(M) \to E_1^{-i,i}(M')$ is epimorphic
for all $i$; the rest of the proof is the same as for Theorem 1.

\subsection{First-order cohomology classes of knots in $\R^n$.}

It is well-known that there are no first-order knot invariants in $\R^3$,
see \cite{v90}. However, the subgroup
$F_{1,\Z_2}^* \subset H^*(\Omega_f \R^3 \sm \Sigma, \Z_2)$
of all $\Z_2$-valued first-order cohomology classes 
is non-trivial: it has exactly
two non-trivial components $F_{1,\Z_2}^1 \simeq F_{1,\Z_2}^2 \simeq \Z_2$.

The generator of the first of them can be defined as 
the linking number with the cycle in $\Sigma$, formed by all maps
$\phi: S^1 \to \R^3$, gluing together some two {\it opposite}
points of $S^1$; the generator of the group $F_{1,\Z_2}^2$ is
just the square of this one.
\medskip

More generally, the following statement holds.
\medskip

{\sc Theorem 2.} {\it For any $n \ge 3$, the subgroup
$F_{1,\Z_2}^* \subset H^*(\Omega_f \R^n \sm \Sigma, \Z_2)$ of 
first-order cohomology classes of the space of knots in $\R^n$ 
contains exactly two non-trivial
components $F_{1,\Z_2}^{n-2} \sim F_{1,\Z_2}^{n-1} \sim \Z_2$.
The generator of the first of them is equal to the linking number 
with the set of maps gluing together some two opposite points
of the circle. The generator of the second  
can be realized by the linking number with 
a similar variety, where these two opposite points are fixed,
say are equal to $0$ and $\pi$, and in the case of odd $n$ is equal
to the Bockshtein of the first generator.

If $n$ is even, then both these cohomology classes give rise to
integer cohomology classes, i.e. $F_{1,\Z}^{n-2} \sim F_{1,\Z}^{n-1} \sim \Z$,
and there are no other non-trivial integer cohomology 
groups $F_{1,\Z}^d$, $d \ne n-2, n-1$.} \quad $\Box$
\medskip

These cocycles in 
$\Omega_f \R^3 \sm \Sigma$ 
are non-trivial already in restriction to
the component of unknots.
Indeed, consider the
standard embedded circle in $\R^3$ and rotate it by all 
angles $\alpha \in [0,2\pi]$ around
any of its diameters. Then we obtain a nontrivial element of 
$H_1(\Omega_f \R^3 \sm \Sigma,\Z_2)$, which takes non-zero value
on the generator of the group $F^1_{1,\Z_2}$.
\medskip

\begin{figure}
\unitlength=1.00mm
\special{em:linewidth 0.4pt}
\linethickness{0.4pt}
\begin{picture}(130,48)
\put(12,23){\oval(24,6)}
\put(18,43){\oval(6,6)[t]}
\put(18,43){\oval(6,6)[bl]}
\put(18,41){\oval(6,2)[br]}
\put(30,43){\oval(18,6)[b]}
\put(30,43){\oval(18,6)[tr]}
\put(30,45){\oval(18,2)[tl]}
\put(52,43){\oval(12,6)[t]}
\put(52,43){\oval(12,6)[bl]}
\put(52,41){\oval(12,2)[br]}
\put(64,43){\oval(12,6)[b]}
\put(64,43){\oval(12,6)[tr]}
\put(64,45){\oval(12,2)[tl]}
\put(85,43){\oval(18,6)[t]}
\put(85,43){\oval(18,6)[bl]}
\put(85,41){\oval(18,2)[br]}
\put(97,43){\oval(6,6)[b]}
\put(97,43){\oval(6,6)[tr]}
\put(97,45){\oval(6,2)[tl]}
\put(107,23){\oval(24,6)}
\put(85,3){\oval(18,6)[b]}
\put(85,3){\oval(18,6)[tl]}
\put(85,5){\oval(18,2)[tr]}
\put(97,3){\oval(6,6)[t]}
\put(97,3){\oval(6,6)[br]}
\put(97,1){\oval(6,2)[bl]}
\put(52,3){\oval(12,6)[b]}
\put(52,3){\oval(12,6)[tl]}
\put(52,5){\oval(12,2)[tr]}
\put(64,3){\oval(12,6)[t]}
\put(64,3){\oval(12,6)[br]}
\put(64,1){\oval(12,2)[bl]}
\put(18,3){\oval(6,6)[b]}
\put(18,3){\oval(6,6)[tl]}
\put(18,5){\oval(6,2)[tr]}
\put(30,3){\oval(18,6)[t]}
\put(30,3){\oval(18,6)[br]}
\put(30,1){\oval(18,2)[bl]}
\put(12,20){\circle*{1.33}}
\put(27,40){\circle*{1.33}}
\put(58,43){\circle*{1.33}}
\put(88,46){\circle*{1.33}}
\put(107,26){\circle*{1.33}}
\put(88,6){\circle*{1.33}}
\put(27,0){\circle*{1.33}}
\end{picture}
\caption{Non-trivial 1-cycle in the space of unknots}
\label{ha}
\end{figure}
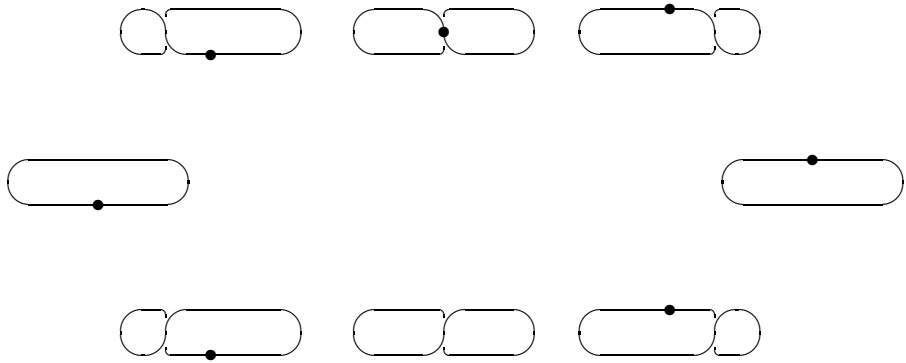

Indeed, let us realize this 1-cycle by the (obviously homotopic
to it) family of unknots shown in fig.~\ref{ha}. Then span it by a disc
in $\Omega_f \R^3$, swept out by the 1-parametric family of segments, 
connecting in the shortest way any two unknots of our family, placed
in this picture one over the other, so that along any such segment
the projection to $\R^2$ remains the same. It is obvious that 
the intersection number of this disc with the above-mentioned subvariety 
in $\Sigma$, generating the group $F_{1,\Z_2}^1$, is non-trivial
(mod 2). 

On the other hand, it is easy to prove that our 1-cycle in the space
of unknots is homotopic there to the cycle, consisting of embeddings
with the one and the same image, which are obtained one from the
another by shifts of the cyclic parameter $\alpha$. 

Consider the space of naturally parametrized great circles in
a sphere $S^2 \subset \R^3$. This space is obviously homeomorphic
to $SO(3) \sim \R P^3$. The restriction on it of our generator of
$F^1_{1,\Z_2}$ coincides with the generator of its $\Z_2$-cohomology
ring, hence also its square (generating $F^2_{1,\Z_2}$) is nontrivial
in restriction to this space.

\section{Construction of the spectral sequence(s)}

\subsection{Resolution spaces.}

Denote by $\Psi$ the space of all unordered pairs of points in
$S^1$ (may be coinciding): $\Psi = S^1 \times S^1/\{\alpha\times \beta =
\beta \times \alpha\}$. It is easy to see that $\Psi$ is 
diffeomorphic to the closed M\"obius band.

Let $\U: \Psi \to \R^\kappa$ be a generic embedding of $\Psi$
into the space of a huge (may be finite) dimension.
For any map $\phi: S^1 \to M$ of the class $\Gamma_U$, consider 
all the points $(\alpha, \beta) \in \Psi$ such that either
$\alpha=\beta$ and $\phi(\alpha)=\phi(\beta)$, or
$\alpha=\beta$ and $\phi'=0$ at the point $\alpha$.
If $\Gamma$ is not very degenerate, then the number of such
points for any 
$\phi \in \Sigma \cap \Gamma_U$ 
is estimated from above by an uniform number (depending on $\Gamma$);
we suppose that the dimension of $\R^\kappa$ is sufficiently 
large with respect to this number.

Consider all images $\U(\alpha,\beta) \in \R^\kappa$ of all such points
for this $\phi$. If $\kappa$ is sufficiently large and the embedding $\U$
is generic, then all these points are vertices of a certain simplex in
$\R^\kappa$; denote this simplex by $\Delta(\phi)$. Define the 
space $\sigma(\Gamma)$ as the subset in $\Gamma_U \times \R^\kappa$
swept out by all simplices of the form $\phi \times \Delta(\phi)$
over all $\phi \in \Sigma \cap \Gamma_U$.
\medskip

{\sc Proposition 5} (cf. \cite{book}). 
{\it The map $\sigma(\Gamma) \to \Sigma \cap \Gamma_U$, defined by the
obvious projection $\Gamma_U \times \R^\kappa \to \Gamma_U$,
is proper, and the induced map
\begin{equation}
\label{two}
\bar H_*(\sigma(\Gamma)) \simeq \bar H_*(\Sigma \cap \Gamma_U)
\end{equation}
is an isomorphism. \quad $\Box$}

\subsection{Filtration and stratification of the resolution set.}

The spaces $\sigma(\Gamma)$ have nice structures which allow 
(in principle) to calculate groups (\ref{two}), namely, the filtration 
$F_1 \subset F_2 \subset \cdots$
(by the ``complexities'' of underlying singularities in $\Sigma$)
and a decomposition of terms $F_i \setminus F_{i-1}$ of this filtration
in correspondence with a certain classification
of these singularities.
These structures have two useful properties:

1) cohomology groups, associated with these objects, are functorial with respect
to embeddings $\Gamma \subset \Gamma'$, see \S \ 2.3 below;

2) these groups converge (in some weak sence) to cohomology
groups of more or less standard topological spaces like the space of continuous 
maps of a given graph to $M$, see \S \ 2.4.
\medskip

The construction of these structures is based on the following 
classification of singular knots.
\medskip

{\sc Definitions} (cf. \cite{v90}, \cite{book}). Let  
$A = \{a_1, \ldots, a_{\# A}\}$ 
be an arbitrary 
finite unordered collection of naturals, 
all whose members $a_l$ are 
not less than 2; let $b$ be an nonnegative integer. Denote by $|A|$ the sum of
all numbers $a_l$. An $A$-{\it configuration}
is any family of $|A|$ pairwise distinct points in $S^1$ partitioned into 
$\# A$ groups of cardinalities $a_1, \ldots, a_{\# A}$ respectively.
An $(A,b)$-{\it configuration} is a pair consisting of an $A$-configuration
and an additional family of $b$ pairwise distinct points in $S^1$ (some of
which can coincide with points of the $A$-configuration).
The map $\phi: S^1 \to M$ {\it respects} the $(A,b)$-configuration
if it sends all points of any of its groups of cardinalities
$a_1, \ldots, a_{\# A}$ into one point in $M$, and $\phi'=0$
at all points of its $b$-part. Two $(A,b)$-configurations are {\it equivalent}
if they can be transformed into one another by an orientation-preserving
homeomorphism of $S^1$.

For instance, the $[i]$-, $\li$- and $i^*$-configurations from \S \ 1.3
are respectively the $(A,b)$-configurations with $A=(2, \ldots,2)$
($i$ twos), $b=0$; $A=(3,2,\ldots, 2)$ ($i-2$ twos), $b=0$, and
$A=(2,\ldots, 2)$ ($i-1$ twos), $b=1$ with the last point different from
$2i-2$ points forming the $A$-part.

The {\it complexity} of an $(A,b)$-configuration 
is the number $|A|-\# A +b$, so that the codimension in $\OM$ of the set
of respecting it maps is equal to $\dim M$ times this number.

For any equivalence class $\J$ of $(A,b)$-configurations, $\delta(\J)$
is the dimension of this class, i.e. the number of geometrically
distinct points in any  configuration $J$ of this class.

An $(A,b,M)$-{\it configuration} is any pair, consisting of an 
$(A,b)$-configuration in $S^1$  and some 
$\# A +b$ points $m_1, \ldots, m_{\# A +b}$ in $M$.
A map $\phi: S^1 \to M$ {\it respects} such a configuration, if it 
sends any group of $a_l$ points, participating in the definition of the 
$A$-part of the configuration, into the point $m_l$, sends any point
$v_j$, participating in the definition of the $b$-part, into $m_{\# A+j}$,
and $\phi'(v_j)=0$ for any such point $v_j$. 
The configuration is {\it acceptable},
if it can be respected by at least one map (this means that
if some of points $v_j$
coincide with the points of the $A$-part, then the corresponding points
$m_*$ also coincide). 
\medskip

{\sc Definition.} 
An affine finite-dimensional subspace $\Gamma$ of the space of smooth
maps $S^1 \to \R^N$ is $(M,d)$-{\it nondegenerate}, if 

a) for any acceptable $(A,b,M)$-configuration of complexity $\le d$, the
set of maps $\phi \in \Gamma_U$, respecting 
this configuration, is a smooth
submanifold in $\Gamma_U$, and differentials of all $n(|A|+2b)$ conditions
distinguishing this manifold (i.e. of conditions $\phi(x_1)= \cdots =
\phi(x_{a_1})=m_1, \ldots, \phi(v_1)=m_{\# A+1}, \phi'(v_1)=0, \ldots $)
are linearly independent at any its point;

b) for any equivalence class $\J$ of $(A,b)$-configurations 
of arbitrary complexity, the codimension
in $\Gamma_U$ of the set of maps, respecting some configurations
of this class, is not less than
$n(|A|-\# A+ b)-\delta(\J)$. 
\medskip

{\sc Proposition 6.} {\it For any $d$, $(M,d)$-nondegenerate 
spaces exist and are 
dense in the space of all affine subspaces of sufficiently large dimension in
$C^\infty(S^1,\R^N)$.} 
\medskip

Indeed, for any natural 
$D$ and any $(A,b,M)$-configuration,
the set of $D$-di\-men\-si\-o\-nal subspaces, not satisfying condition
a) at this configuration, is a subvariety in the space of
all $D$-dimensional subspaces in $\Omega_f\R^N$. 
The codimension of this subvariety grows to infinity together
with $D$, in particular for large $D$ becomes greater than the 
dimension of the space of all $(A,b,M)$-configurations with
given $A$ and $b$. Condition
b) follows from the Thom transversality theorem, cf. \cite{book}.
\quad $\Box$ \medskip

{\sc Definition.} For any $(A,b)$-configuration $J$, the simplex 
$\Delta(J) \subset \R^\kappa$ is defined as the simplex $\Delta(\phi)$
for any generic $\phi \in \OM$ respecting $J$
(i.e. having no extra singularities).
\medskip

The number of vertices of this simplex is
equal to $\sum_{j=1}^{\# A} {a_j \choose 2}+b$.
\medskip

{\sc Definition.} Given any equivalence class $\J$
of $(A,b)$-configurations, the corresponding
$\J$-{\it block} $B(\J,\Gamma)$ is the union of all points 
$(\phi,\zeta) \in \sigma(\Gamma) \subset \Gamma_U \times \R^\kappa$
such that for some $(A,b)$-configuration $J \in \J$

a) $\phi$ respects $J$, and

b) $\zeta$ belongs to the simplex $\Delta(J)$.
\medskip

The term $F_i$ of the {\it main filtration} of $\sigma(\Gamma)$
is defined as the union of all $\J$-blocks over all $\J$ of
complexity $\le i$.

By definition, $B(\J,\Gamma)$ consists of simplices 
$\zeta \times \Delta(J) \sim \Delta(J)$,
$J \in \J$. Any such simplex lies in $F_i$ (where $i$ is the complexity
of $\J$), but some of its points belong to $F_{i-1}$. These points
constitute several faces of $\Delta(J)$. Namely, any face of $\Delta(J)$
is characterized by the collection of its vertices, i.e. by a collection of
$\# A$ graphs with $a_1, \ldots, a_{\# A}$ vertices respectively,
and a choice of some of $b$ ``singular'' points of the configuration.
The faces lying in $F_{i-1}$ are exactly those 
that either one of corresponding graphs is not connected, or
at least one of $b$ points is missed in this choice, cf. \cite{v90}.
\medskip

{\sc Proposition 7} (see e.g. \cite{book}). {\it For any $(A,b)$-configuration
$J$ of complexity $i$, the group
$\bar H_*(\Delta (J) \setminus F_{i-1})$ is trivial in all
dimensions other than $\sum_{j=1}^{\# A}(a_j-1)+b-1 \equiv i-1,$
and in dimension $i-1$ it is isomorphic to
\begin{equation}
\label{tens}
\otimes_{j=1}^{\# A}\Z^{(a_j-1)!} 
\end{equation}
(in particular to $\Z$ if all
$a_j$ are equal to 2).}
\medskip

For $\J$ of complexity $i$ denote by $\tilde B(\J,\Gamma)$
the ``pure part'' $B(\J,\Gamma)\setminus F_{i-1}$ of the
$\J$-block $B(\J,\Gamma)$. By the construction, it is 
the space of a fiber bundle, whose base is the 
space of all pairs of the form \{an $(A,b)$-configuration $J \in \J$, a map 
$\phi \in \Gamma_U$ respecting $J$\}, and the fiber over such a point 
is the set of interior points of 
faces of the simplex $\Delta(J)$, not belonging to $F_{i-1}$,
so that the Borel--Moore homology group of the fiber is
described by Proposition 7.

Denote by $\beta(\J,\Gamma)$ the base of this fiber bundle,
and by $\beta(\J)$ the space of similar pairs $\{J,\phi\}$
over all $J \in \J$ and all $\phi \in \OM$ respecting $J$
(and not only over such $\phi \in \Gamma_U$).
\medskip

{\sc Proposition 8.} {\it 1.  
For any equivalence class $\J$ of $(A,b)$-configurations 
of complexity $d$
and any homology class $\xi \in H_*(\beta(\J))$ 
(with coefficients in any local system of groups) there exists 
a $(M,d)$-nondegenerate space $\Gamma$ such that $\xi$
can be realized by a cycle belonging to $\beta(\J, \Gamma)$.

2. If two such cycles in $\beta(\J, \Gamma)$ define the same homology
class in $\beta(\J)$, then
there exists a 
$(M,d)$-nondegenerate
space $\Gamma'$, containing $\Gamma$, such that these cycles are homological
already in $\beta(\J, \Gamma')$.}
\medskip

{\it Proof}. This follows from the Weierstrass approximation theorem:
it is sufficient to take weakly moved spaces of maps $S^1 \to \R^N$
given by trigonometric polynomials of sufficiently large degrees.
\quad $\Box$ \medskip

{\sc Definition.} A sequence of finitedimensional affine
subspaces $\Gamma^1 \subset 
\Gamma^2 \subset \cdots $ in $ \Omega_f \R^N$ is {\it exhausting} if

a) for any $d$ almost all its terms $\Gamma^j$ (i.e. all 
except may be for finitely many)  
are $(M,d)$-nondegenerate; 

b) for any $\J$ and any
class $\xi \in H^*(\beta(\J))$, the  condition 1) of the
previous proposition is satisfied for almost all $\Gamma^j$;  

c) for any term $\Gamma^j$ of this sequence and any 
two cycles $\xi, \zeta \in H_*(\beta(\J,\Gamma^j)),$
defining the same element of $H_*(\beta(\J))$,
these cycles are homological in almost all spaces $\beta(\J,\Gamma^k),$
$k \ge j$.
\medskip

Proposition 8 implies that such sequences exist.

\subsection{Stabilization of spectral sequences.}

For any subspace $\Gamma_U \subset \Omega_f \R^N,$ consider the homological
spectral sequence $E^r_{p,q}(\Gamma)$, converging to the group
$\bar H_*(\sigma(\Gamma)) \equiv \bar H_*(\Sigma \cap \Gamma_U)$
and defined by the main filtration of $\sigma(\Gamma)$, described in the
previous subsection. By definition, 
$E^1_{p,q}(\Gamma) \simeq \bar H_{p+q}(F_p(\sigma(\Gamma)) \setminus
F_{p-1}(\sigma(\Gamma)).$ 

Using the formal inversion (\ref{inver}), we convert it to the 
{\it cohomological} spectral sequence
$E_r^{p,q}(\Gamma) \to H^{p+q+1}(\Gamma_U, \Gamma_U \setminus \Sigma).$

For any $d$, denote by $_dE_r^{p,q}(\Gamma)$ the truncated 
spectral sequence, obtained from the previous one by replacing
by 0 all terms $_dE_1^{p,q}$ with $p< -d$. It converges to the 
Borel--Moore homology group of $F_d(\sigma(\Gamma))$:
$_dE_r^{p,q}(\Gamma) \to \bar H_{\dim \Gamma-p-q-1}(F_d(\sigma(\Gamma)))$.

Let $\Gamma \subset \Gamma'$ be two $(M,d)$-nondegenerate 
subspaces in $\Omega_f\R^N$.
Then there is a natural homomorphism
\begin{equation}
\label{stab}
 _dE_r^{p,q}(\Gamma') \to _dE_r^{p,q}(\Gamma).
\end{equation}
Indeed, by the definition of $(M,d)$-nondegeneracy 
$F_d(\sigma(\Gamma))$ admits in
$F_d(\sigma(\Gamma'))$ a tubular neighborhood, 
homeomorphic to the direct product of $F_d(\sigma(\Gamma))$ and
an open $(\dim \Gamma' - \dim \Gamma)$-dimensional disc,
in such a way  that this homeomorphism preserves natural filtrations
of both spaces.

The homomorphism (\ref{stab}) is defined as the 
composition of the restriction on this
neighborhood and the K\"unneth formula in it.
\medskip

The {\it stable spectral sequence}  
$ E_r^{p,q} $ is defined by 
$E_r^{p,q} \equiv$
lim ind $_dE_r^{p,q}(\Gamma^j)$ over such homomorphisms for any
exhausting sequence of approximating spaces $\Gamma^j$.

It is easy to see that it does not depend on the choice of this sequence of
spaces and converges to some subgroup in $H^*(\OM,\OM \setminus \Sigma)$.
\medskip

{\sc Proposition 9.} {\it The support of the stable sequence $E_r^{p,q}$
(i.e. the set of such $p,q$ that $E_1^{p,q} \ne 0$) belongs
to the wedge from fig. 1: $p <0, q+(n-2)p \ge 0$.}
Moreover, the same is true for any non-stable spectral sequence 
$E_r^{p,q}(\Gamma)$ with any $(M,d)$-nondegenerate $\Gamma$.
\medskip

{\it Proof.} By Proposition 7, for any equivalence class $\J$ of 
$(A,b)$-configurations of
complexity $i$ and any element $\Gamma^j$ of our exhausting sequence,
the contribution of the block 
$\tilde B(\J,\Gamma^j) \subset \sigma(\Gamma^j)$ into
the group $E_1^{-i,q} \simeq 
\bar H_{\dim \Gamma -q+i-1}(F_d \setminus
F_{d-1})$ can be nontrivial only if
$\dim \Gamma -q+i-1 \le \dim \beta(\J, \Gamma)+
\sum_{k=1}^{\# A}(a_k-1)+b-1 \le \dim \Gamma -ni+2i+i-1$;
the summand $2i$ in the last expression is the upper estimate 
for the dimension $\delta(J)$ of the space of $(A,b)$-configurations
of the class $\J$. 
This implies the first statement of the proposition. 
The second statement follows
in the same way from condition b) of the definition of 
$(M,d)$-nondegenerate spaces.
\quad $\Box$ \medskip

In the next subsection we show that such stable 
sequences are not too wild and huge.

\subsection{On the calculation of the term $E_1^{p,q}$ of the stable 
spectral sequence.}

Let us fix a certain natural $i$.

By the definition of spectral sequences,
\begin{equation}
\label{e1}
E_1^{-i,q}(\Gamma) \simeq \bar H_{\dim \Gamma+i-q-1}(F_i\setminus F_{i-1}),
\end{equation}
where $F_i \equiv F_i(\sigma(\Gamma))$. The space 
$ F_i \setminus F_{i-1}$ splits into pure $\J$-blocks $\tilde B(\J,\Gamma)$
with $\J$ of complexity $i$. As in \cite{v90}, \cite{book},
we introduce the {\it auxiliary 
filtration} in it, defining its term $\Phi_\alpha$ as the union
of all blocks $\tilde B(\J,\Gamma)$ such that $\delta(\J) \le \alpha$, i.e. 
the configuration $\J$ consists
of $\le \alpha$ geometrically distinct points.
Let $\E^\rho_{\mu,\nu}(\Gamma,i)$ be the spectral sequence, converging
to the group $\bar H_*(F_i \setminus F_{i-1})$ and 
generated by this filtration.
Its term  $\E^1_{\mu,\nu}$ is the direct sum of groups
$\bar H_{\mu+\nu}(\tilde B(\J,\Gamma))$ over all $\J$ of complexity 
$i$ and $\delta(\J)=\mu$.

The stabilization of this spectral sequence 
over growing $\Gamma$ is defined in the same way as 
for the main spectral sequence and allows us to
define the stable cohomological spectral sequence 
$\E_\rho^{a,b}(i) \equiv 
\lim $ ind $\E^\rho_{-a, \dim \Gamma^k-b-1}(\Gamma,i)$
over any exhausting sequence $\{\Gamma^k\}$.
\medskip

{\sc Proposition 10.} {\it 
The stable auxiliary spectral sequence converges to the
term $E_1$ of the stable main spectral sequence.
Namely, its group $\oplus_{a+b=t}\E^{a,b}_\infty(i)$
is adjoined to the group $E^{-i,i+t}_1$.
} \quad $\Box$
\medskip

On the other hand, the stable member $\E_1^{a,b}(i)$
of this sequence can be expressed in terms of cohomology
groups of spaces $\beta(\J)$. Indeed, 

a) this term splits into the direct sum of certain homology groups,
associated with all $\J$ with  complexity 
$i$ and $\delta(\J)=-a$;

b) namely, for any such $\J$ the corresponding summand is the stabilization
of groups $\bar H_{\dim \Gamma^k -a-b -1}(B(\J,\Gamma^k))$;

c) since $B(\J, \Gamma^k)$ is a fiber bundle with base
$\beta(\J,\Gamma^k)$ and fiber described in
Proposition 7,
these stable homology groups are isomorphic to the stabilization of
$(\dim \Gamma^k -a-b-i)$-dimensional Borel--Moore homology groups
of $\beta(\J,\Gamma^k)$ with coefficients in a certain local system
$\Xi$ with fiber (\ref{tens});

d) by the definition of $(M,d)$-nondegeneracy, all spaces
$\beta(\J,\Gamma^k)$ with sufficiently large $k$ are smooth
$(\dim \Gamma^k-a-n\cdot i)$-dimensional 
manifolds, therefore  by the Poincar\'e duality theorem
previous homology groups are isomorphic to groups
$$H^{b-(n-1)i}(\beta(\J,\Gamma^k), \Xi^* \otimes Or(\J,k)),$$
where $Or(\J,k)$ is the orientation sheaf of the manifold 
$\beta(\J,\Gamma^k)$;

e) these sheaves $Or(\J,k)$ stabilize to a common sheaf $Or(\J)$
on $\beta(\J)$ (i.e., they are preserved by all inclusions
$\beta(\J,\Gamma^k) \to \beta(\J,\Gamma^l)$ with sufficiently
large $k<l$);

f) and finally we obtain from Proposition 8 that 
the summand in $\E_1^{a,b}(i)$, corresponding to the
class $\J,$ is isomorphic to
\begin{equation}
\label{poin}
H^{b-(n-1)i}(\beta(\J),\Xi^* \otimes Or(\J)).
\end{equation}

Certainly, it splits into the direct product of similar groups over all
path-com\-po\-nents of $\beta(\J)$, i.e., over the $\J_M$-{\it routes},
cf. \S \ 1.3.
\medskip

{\sc Example} (cf. \S \ 1). Let $n=3$. To calculate the groups $E_r^{-i,i}$,
we need to consider only the groups $\E_\infty^{a,b}$ with $a+b=0$
or $a+b=1$. First suppose that $a+b=0$. Since $a \ge -2i$, we have
$b \le 2i$, thus the group (\ref{poin}) can be nontrivial only if
$b=2i$, i.e., $a=-2i, \mu=2i,$ and $\J$ is a class of $[i]$-configurations.
It is easy to calculate that for such $\J$ the sheaf $\Xi^* \otimes Or(\J)$
is isomorphic to $\Z$ if $M$ is orientable, and coincides with the sheaf
$s\Z$ from formula (\ref{nonor}) otherwise.

Thus the unique group $\E_1^{a,b}$ with $a+b=0$ is the group
$\E_1^{-2i,2i} \equiv H^0(\beta(\J),s\Z)$, formally 
generated by $s$-orientable $[i]_M$-routes.

Similarly, the only two non-trivial groups $\E_1^{a,b}$ with $a+b=1$
are the group $\E_1^{-2i,2i+1}$, which is just the group
(\ref{nonor}) (or (\ref{circ}) if $M$ is orientable) and the group
$\E_1^{-2i+1,2i}$, which is the direct sum of certain 0-dimensional
homology groups of $\li_M$- or $i^*$-routes. In the case of 
$\li_M$-routes, the coefficient sheaf $\Xi^*$ is locally isomorphic
to $\Z^2$ and consists of $\Z$-linear combinations of corresponding
$\ti_M$-routes factorized through the diagonal, consisting of such
combinations with coinciding coefficients.

\medskip
{\sc Concluding remark.} All the above considerations 
can be word-for-word carried out to the case 
of $d$-component links in $M$ with any fixed $d$:
in the basic construction instead of the space of maps
of one circle to $M$ (or to $\R^N$) we need only to consider
the space of similar maps of the disjoint union of
$d$ circles, cf. \cite{stan}

\bigskip

\end{document}